\documentclass{imsart}
\usepackage{amsmath,amssymb,amsthm,booktabs}
\usepackage{xr}
\usepackage{verbatim}
\usepackage{graphicx}
\usepackage[OT1]{fontenc}
\usepackage[colorlinks,citecolor=blue,urlcolor=blue]{hyperref}
\usepackage{txfonts}
\usepackage{indentfirst}
\usepackage{multirow}
\usepackage{float,subfig}
\usepackage{dsfont}
\usepackage{bm}
\usepackage{xcolor}
\usepackage{graphicx}
\usepackage[round]{natbib}
\usepackage{enumitem}
\startlocaldefs
\numberwithin{equation}{section} \theoremstyle{plain}
\newtheorem{theorem}{Theorem}[section]
\newtheorem{lemma}{Lemma}[section]
\newtheorem{corollary}{Corollary}[section]

\newtheorem{definition}{Definition}
\newtheorem{remark}{Remark}[section]
\endlocaldefs

\DeclareMathOperator{\tr}{Tr}
\DeclareMathOperator{\cov}{Cov}

\makeatletter
\def\hlinewd#1{%
	\noalign{\ifnum0=`}\fi\hrule \@height #1 \futurelet
	\reserved@a\@xhline}
\makeatother

\begin{document}
	
	
	\newcommand\E{\mathbb{E}}
	\newcommand\lb{\left(}
	\newcommand\rb{\right)}
	\newcommand{\lv}{\left\lvert}
	\newcommand{\rv}{\right\rvert}
	\newcommand\z{\mathbf{z}}
	\newcommand{\veps}{\varepsilon}
	\newcommand\norm[1]{\left\lVert#1\right\rVert}
	\newcommand\x{\mathbf{x}}
	\newcommand{\y}{\mathbf{y}}
	\newcommand\vb{\mathbf{\overline{v}}}
	\newcommand\uu{\mathbf{u}}
	\newcommand{\vv}{\mathbf{v}}
	\newcommand{\e}{\mathbf{e}}
	\newcommand\T{\mathsf{T}}
	
	\renewcommand\a{\mathbf{a}}
	\renewcommand{\d}[1]{\ensuremath{\operatorname{d}\!{#1}}}
	
	\allowdisplaybreaks \allowdisplaybreaks[1]
	
	\newcommand\gai[1]{{\color{red}#1}}
	\def\bl{\color{blue}}

	\def\G{\mathbf{G}}
	\def\U{\mathbf{U}}
	\def\K{\mathbf{K}}
	\def\A{\mathbf{A}}
	\def\B{\mathbf{B}}
	\def\C{\mathbf{C}}
	\def\H{\mathbf{H}}
	\def\W{\mathbf{W}}
	\def\D{\mathbf{D}}
	\def\X{\mathbf{X}}
	\def\Y{\mathbf{Y}}
	\def\I{\mathbf{I}}
	\def\N{\mathbf{N}}
	\def\Z{\mathbf{Z}}
	\def\TT{\mathbf{T}}
	\def\WW{\mathcal{W}}
	\def\S{\mathbf{S}}
	\def\R{\mathbf{R}}

	\def\Thta{\mathbf{\Theta}}
	\def\th{\mathbf{\theta}}

	\def\GA{\mathbf{\Gamma}}
	\def\SIG{\mathbf{\Sigma}_{\x}}
	
	\def\SI{\mathbf{\Sigma}}
	\def\diag{\text{diag}}
	\def\s{\underline{s}}
	\def\m{\underline{m}}
	
	\renewcommand{\baselinestretch}{1.2}
	
	\begin{frontmatter}
		
		\title{Central Limit Theorem for Linear Spectral Statistics of
			Large Dimensional Kendall's Rank Correlation Matrices and its Applications}
		\runtitle{CLT for LSS of Rank Correlation Matrix}

		\begin{aug}
			\author{\fnms{Zeng} \snm{Li}\thanksref{m1}\ead[label=e1]{zxl278@psu.edu}},
			\author{\fnms{Qinwen} \snm{Wang}\thanksref{m2}\ead[label=e2]{wqw@fudan.edu.cn}}
			\and
			\author{\fnms{Runze} \snm{Li}\thanksref{m3}
				\ead[label=e3]{rzli@psu.edu}}
			\thankstext{m1}{Zeng Li's research is supported by NIDA, NIH grants
				P50 DA039838.}
			\thankstext{m2}{Wang's research is partially supported
				by NSFC (No.\ 11801085) and Shanghai Sailing Program (No.\ 18YF1401500).
				Wang is the corresponding author.}
			\thankstext{m3}{Li's research is supported by NIH grants
				P50 DA039838 and R01 ES019672 and  a NSF grant DMS 1820702.}
			
			\runauthor{Z. Li, Q. Wang and R. Li}
			
			\affiliation{Pennsylvania State University\thanksmark{m1}\thanksmark{m3} and
				Fudan University\thanksmark{m2}}
			
			\address{Zeng Li\\
				Department of Statistics\\
				Pennsylvania State University\\
				University Park, PA, 16802\\
				\printead{e1}}

			\address{Qinwen Wang\\
				School of Data Science\\
				Fudan University\\
				Shanghai, China\\
				\printead{e2}}
			
			\address{Runze Li\\
				Department of Statistics\\
				Pennsylvania State University\\
				University Park, PA, 16802\\
				\printead{e3}}
		\end{aug}

\begin{abstract}
This paper is concerned with the limiting spectral behaviors of large
dimensional  Kendall's  rank correlation matrices generated by samples 
with independent and continuous components. We do not require the components  to be identically distributed, and do not need any moment conditions, which is very different from the assumptions imposed in the literature of random
matrix theory. The statistical setting in this paper covers a wide range of
highly skewed and heavy-tailed distributions. We establish the central limit
theorem (CLT)  for the linear spectral statistics  of the  Kendall's rank
correlation matrices under the Marchenko-Pastur asymptotic regime, in which  the dimension 
diverges to infinity proportionally with
the sample size.
We further  propose three nonparametric procedures for high dimensional independent test and their limiting null distributions are derived by 
implementing this CLT. Our numerical comparisons demonstrate the robustness and superiority of
our proposed test statistics under various mixed and heavy-tailed cases.
    \end{abstract}
    
 \begin{keyword}[class=MSC]
    \kwd[Primary ]{62H10}
    \kwd[; secondary ]{62H15}
  \end{keyword}

  \begin{keyword}
    \kwd{Kendall's rank correlation matrices}
    \kwd{linear spectral statistics}
    \kwd{central limit theorem}
    \kwd{random matrix theory}
    \kwd{high dimensional independent test}
  \end{keyword}

\end{frontmatter}

\section{Introduction}\label{sec:intro}

Most well-established statistics in classical multivariate analysis can be
presented as linear functionals of eigenvalues of sample covariance or
correlation matrix models. Such linear functionals of eigenvalues are termed as
{\em linear spectral statistics (LSS)} in the literature of {\em Random Matrix Theory (RMT)}.
Studies of asymptotic properties of the LSS are particularly important in
multivariate analysis of variance, multivariate linear models, canonical
correlation analysis and factor analysis, etc. Analysis of high dimensional data
calls for new theoretical framework  since classical multivariate
analysis theory may become invalid under high-dimensional regime. While 
RMT, which studies asymptotic behaviors of eigenvalues  of large random matrices with certain structures, serves as an effective
tool to deal with high dimensional problems, particuarly under the {\em Marchenko-Pastur asymptotic regime}, where the dimension proportionally
diverges to infinity with
the sample size. In general, the study of RMT can be divided
into two categories: bulk spectrum and edge behaviors. The bulk spectrum
describes global eigenvalue behaviors including convergence of {\em empirical
spectral distribution (ESD)}, fluctuation of LSS and spectrum separation, etc; while
edge behaviors mainly concern the convergence and fluctuation of extreme eigenvalues.


The study of the fluctuation  of  LSS for different types of random
matrix models  has received considerable attention in past decades and has
broad applications in various fields such as wireless communications and
finance, etc, see recent monographs  \citep{BSbook,
couillet2011random, yao2015sample} and surveys \citep{tulino2004random, johnstone2006high,
paul2014random}. The original result on CLT for LSS can be traced back to
\citet{jonsson1982some}, where the CLT for polynomial functions  of  a sequence
of  Wishart matrices was established. \citet{BS04} extended it to analytic
functions of eigenvalues of large sample covariance matrices. Theoretical
results for other popular random matrices have been derived {later on, which} include
\citet{zheng12} for the Fisher matrices; \citet{yang15} for canonical
correlation matrices and \citet{gao2017high} for the Pearson correlation
matrices. Other than Pearson-type models,
  \cite{diaconis2001linear} adopted the moment method to establish the CLT for  Haar matrices.   \citet{BY05} obtained the CLT for LSS of Wigner matrices using the Stieltjes transform method.

Although the result in \citet{BS04}  is universal in the sense that the CLT
does not depend on any particular  distribution assumption of the data, it
requires the same kurtosis as Gaussian distribution. The moment restrictions
have been relaxed in subsequent developments like \citet{pan2008central},
\citet{zheng2015substitution}, finite fourth or even higher order moments are
still generally required in Pearson-type models. Such constraints significantly
limit the applicability for real data analysis, but certainly have stimulated
the investigation of  nonparametric methods which are free of moment restrictions. For starters, \citet{B15}
established the asymptotic normality of polynomial functions of the spectrum of
Spearman's rank correlation matrices. \citet{LSD17} obtained the limiting spectral distribution of Kendall's rank
correlation matrices, which appeared to be a variation of the
Marchenko-Pastur law. \citet{BZG18} further proved the Tracy-Widom limit for largest eigenvalues of Kendall's rank correlation matrices. 

In this paper, we focus on LSS of  Kendall's rank correlation matrices $\K_n$ generated by samples  $\x_1, \cdots, \x_n$ from a $p$-dimensional random vector $\x$.
The primary goal of this paper is to  establish the CLT for general LSS of
$\K_n$ under the Marchenko-Pastur asymptotic regime when $\x$ consists of $p$
independent components.
From technical point of view,
the LSS of $\K_n$ is actually the type of U-statistics and to the best of our
knowledge, there are no such  CLTs for LSS related to U-statistics have been established so far in the
literature. Moreover, the structure of $\K_n$ is very different from most of the well studied random
matrix models in the literature. Although $\K_n$ can be
written as Wishart-type product of data matrix $\Thta$, i.e.,
$\K_n=\Thta\Thta^{\T}$ (see more details in \eqref{eq:Kn}), the elements
within each row of $\Thta$ is nonlinearly correlated and such correlation
structure cannot be incorporated into any random matrix models studied so far. In fact, the most commonly used model assumption is the {\em independent
component model (ICM)} assumption: there exist constant vector $\bm{\mu}$ and
non-negative definite matrix $\SI_p$ such that $\x=\bm{\mu} + \SI^{1/2}_p{\bf z},$
where all elements within
${\bf z}$ are independent and identically distributed and have zero mean and unit variance 
\citep{BS04, pan2008central, zheng12}.  Thus, techniques  for the ICM cannot
be directly adopted to accommodate the nonlinear correlation structure in
$\Thta$.
In order to decouple such nonlinear correlated structure of our
target matrix $\K_n$, we develop new technical tools
for each row $\vv^{\T}_k$ of
$\Thta$. Here the tricky thing is that, although $\vv_k^{\T}$ can be decomposed
into two uncorrelated parts using Hoeffding's decomposition, the two components
are still not independent. Due to this kind of subtle dependent relation, it
requires much effort to derive the expectation of product of two quadratic
forms like the term $\E \vv_k^\T \A \vv_k \vv^\T_k \B \vv_k$. In fact, derivation of
explicit formula for such expectation is one of the key steps for deriving the CLT
for LSS of $\K_n$ and the result is established in Lemma~\ref{lem:quadratics}.
Moreover, parallel to the population covariance matrix $\SI_p$ in \citet{BS04},
$M\E(\vv_k\vv_k^{\T})$ is a $M \times M$ ($M=n(n-1)/2$) matrix with $n-1$
eigenvalues equal to $(n+1)/3$ and all the remaining ones being $1/3$. This fact
suggests that the spectral norm of $M\E(\vv_k\vv_k^{\T})$ is unbounded  while
its empirical spectral distribution degenerates to a point mass at $1/3$
asymptotically. Thus, such kind of spectral property differs from
the common settings in the study of  high dimensional sample covariance matrices where its population version $\SI_p$ is of full rank
with bounded spectral norm and its empirical spectral distribution converges to
certain proper c.d.f. asymptotically. Due to this reason,  inequalities related to
higher order moments of $\vv_k \A \vv_k^\T$ and its centralized versions need
to be re-considered and re-established in order to derive the CLT for LSS of
$\K_n$, which are obtained in Lemma~\ref{lem:generalvkbound} and
Lemma~\ref{lem:vkbound}.

To demonstrate the potential of the newly established CLT, we study high dimensional independent tests without any moment
conditions.
Under Gaussian assumption, such tests of independence have been studied in the
low dimensional  framework \citep{Muirhead82, Anderson84},
where  the {\em likelihood ratio test statistic (LRT)} was shown to converge to 
chi-square distribution under the null hypothesis. However, the well established large sample
theory for the LRT becomes invalid when the dimension grows
comparable or even larger than the sample size. To cope with high
dimensionality,  various methods have been proposed in the literature,
including the modifications of LRT \citep{jiang13, jiang2015likelihood};
Frobenius-norm-type statistics based on sample correlation matrices
\citep{gao2017high}, sample canonical correlation coefficients \citep{yang15}
and maximum-norm-type statistics  \citep{jiang04, zhou07, cai2011limiting}.
These tests are in general infeasible for heavy-tailed distributions since
these approaches are based on Pearson-type sample covariance, correlation or
canonical correlation matrices which all require strong moment conditions.

On the other hand, nonparametric approaches can get rid of such moment restrictions and robustness are thus achieved by 
replacing the Pearson-type correlations with rank-based versions. Recent works on this topic include \citet{B15} for tests based on
the polynomial functions of the spectrum of Spearman's rank correlation matrices,
 \cite{han17} for the maximum type statistics  including Kendall's tau and
Spearman's rho, and \citet{Leung18} for a class of nonparametric U-statistics.
In this paper, we propose three test statistics for high dimensional independent test.  The first two  are based on the second and fourth spectral moment of $\K_n$ while the third is based on the {\em entropy loss (EL)} between the Kendall's rank correlation matrix $\K_n$ and its population counterpart as motivated by \citet{james1961estimation},~\citet{muirhead2009aspects},~\citet{zhenghypothesis}. All the three statistics can be represented in certain forms of LSS of $\K_n$, thus their  limiting null distributions can be fully derived using our newly
established CLT for general LSS of $\K_n$. 
Our numerical studies
demonstrate robustness of the three tests
for  mixed and heavy-tailed distributions. Their performances are very satisfactory for distinguishing various settings of dependent alternatives including both linear and non-linear cases.

Throughout this paper, we use bold Roman capital letters to represent matrices, e.g. $\A$. $\tr(\A)$  denotes the trace of matrix $\A$.  Scalars are in lowercase letters.  Vectors are bold letters in lowercase like $\mathbf{v}$. $\mathbb{N}$, $\mathbb{R}$, and  $\mathbb{C}$ represent  the sets of natural, real  and complex numbers. $\E(\cdot)$ means taking expectation and $\Im(\cdot)$ means taking imaginary part of complex numbers.  $\mathds{1}(\cdot)$ stands for indicator function and $\T$ for transpose of vectors or matrices.

The rest of the paper is organized as follows. Section~\ref{kenmatrix}
first introduces some preliminary results and notations on high dimensional Kendall's
rank correlation matrices and then establishes our main result on
the CLT for LSS of $\K_n$. Section~\ref{appli} develops three statistics for  high dimensional
independent test, with their limiting null distributions explicitly derived.  Simulation experiments are conducted to compare the finite sample
performance of our test statistics with other existing ones under various model settings and alternatives.
Section~\ref{mainre} contains the proof of our main result on CLT for LSS of $\K_n$. Technical lemmas
and related proofs are relegated to the supplementary file. 
\section{CLT for  LSS of High-dimensional Kendall's rank correlation matrices}\label{kenmatrix}

Suppose  $\x_i=\big( x_{1i}, x_{2i},\cdots,x_{pi}\big)^{\T}, ~i=1, \cdots,
n$ is a random sample from a $p$ dimensional random vector $\x$ with independent and continuous components. Denote $\X_{n}=\lb
\x_1,\x_2,\cdots,\x_n\rb$ as the data matrix. The Kendall's rank
correlation matrix of $\X_n$ is defined to be  a $p\times p$ matrix
$\K_n=(\tau_{k \ell})$, whose $(k, \ell)$-th entry  is the empirical Kendall's 
rank correlation coefficient between the $k$-th and $\ell$-th components of $\x$ (the $k$-th and $\ell$-th row of $\X_{n}$), i.e.
\begin{align}\label{eq: taukl}
\tau_{k \ell}=\frac{1}{{n\choose 2}}\sum_{1\leq i<j\leq
n}\text{sign}(x_{ki}-x_{kj})\cdot \text{sign}(x_{\ell i}-x_{\ell j})~\quad
1\leq k, \ell \leq p.
\end{align}
For any analytic function $f(x)$,  a general LSS of $\K_n$ is defined as $\sum_{i=1}^pf(\lambda_i)/p$, where $\{ \lambda_i,~1\leq i\leq p\} $ are eigenvalues of $\K_n$.
 In this section, we are interested in the 
CLT for such LSS of $\K_n$ under the Marchenko-Pastur asymptotic regime:
$ p \to \infty, n\to \infty, c_n:=p/n \to c\in(0,\infty).
$
This asymptotic regime is denoted by $(p, n)\to \infty$ for short throughout
the paper.

\subsection{Preliminary results on the limiting spectrum of $\K_n$}
To study the limiting properties regarding to the spectrum  of certain large dimensional random matrix models, the very first step is to consider  its empirical spectral distribution. Suppose the random matrix model that we are interested in is a $n\times n$  Hermitian (or symmetric) matrix $\A_n$, we denote its  $n$ real eigenvalues as $\lambda_1\geq  \cdots \geq \lambda_n$  in descending order. The {\em empirical spectral distribution} ({\em ESD}) of $\A_n$ is referred as a random measure $F^{\A_n}$ such that
\begin{align}\label{anesd}
F^{\A_n}(x)=\frac{1}{n}\sum_{i=1}^n\delta_{\lambda_i}~,
\end{align}
where $\delta_{\lambda_i}$ is  the Dirac mass at point $\lambda_i$.
If there exists a non-random proper c.d.f $F(x)$ such that as $n\to \infty$, with probability 1,
$F^{\A_n}(x)\xrightarrow{d} F(x)$, then $F(x)$ is called the
{\em limiting spectral distribution} ({\em LSD}) of $\A_n$.


The study of the limiting spectrum of $\K_n$ from the perspective of RMT is few and far between. The only two relevant references are
\citet{LSD17} and \cite{BZG18}.  \citet{LSD17} studied the LSD of $\K_n$ and
showed that $F^{\K_n}$ converges in probability to $F_c:=2Y/3+1/3$ as $(p,n)\rightarrow \infty$
where $Y$ follows the standard Marchenko-Pastur law with parameter
$c$,  see \citet{marchenko1967distribution}. This cumulative distribution function $F_c$ has an explicit form of density
 $d F_c$ given by
\begin{equation}\label{kendensity}
d F_c(x)=\frac{9}{4\pi c(3x-1)}\sqrt{\lb d_{c,+} -x\rb \lb x-d_{c,-}\rb}+\Big(
1-\frac{1}{c}\Big)~ \delta_{\frac{1}{3}}\mathds{1}_{\{c>1\}},\quad~d_{c,-}\leq
x\leq d_{c,+},
\end{equation}
where
\begin{equation}\label{eq:dcs}
d_{c,-}= \frac{1}{3}+\frac{2}{3}\lb 1-\sqrt{c}\rb^2,\quad d_{c,+}=\frac
13+\frac{2}{3}\lb 1+\sqrt{c}\rb^2.
\end{equation}
A very useful tool for deriving such LSD is its corresponding {\em
Stieltjes transform}, for any c.d.f $F(x)$, 
  defined  to be
\[m_F(z)=\int \frac{1}{x-z}\d F(x),\quad z\in \mathbb{C}^{+},\]
where $\mathbb{C}^{+}$ denotes the upper complex plane. The Stietjes transform
$m_{F_c}(z)$ of $F_c$ has been proven to be the unique solution of the following
equation
\begin{equation}\label{kenmz}
\frac{2}{3}c\lb z-\frac{1}{3}\rb m^2_{{F_c}}(z)+\lb z-1+\frac{2}{3}c\rb
m_{F_c}(z)+1=0~
\end{equation}
such that $\Im (z)\cdot \Im(m_{{F_c}}(z))>0$. It can also be expressed explicitly as a function of the limiting dimension-to-sample size ratio $c$, i.e.
\begin{align*}
m_{{F_c}}(z)=\frac{1-\frac 23c-z+\sqrt{\big(z-1-\frac 23
c\big)^2-\frac{16}{9}c}}{\frac 43 c (z-\frac 13)}~.
\end{align*}
Apart from such global behavior, recently, \cite{BZG18} proved the Tracy-Widom
law for the largest eigenvalue of $\K_n$ under the same
assumptions,  which turns out to be the first Tracy-Widom law for a
high-dimensional U-statistic.

\subsection{Hoeffding decomposition}
For the $p\times n$ data matrix $\X_n=(\x_1, \cdots, \x_n)$ with $\x_i=\big( x_{1i}, x_{2i},\cdots,x_{pi}\big)^{\T},  i=1, \cdots,
 n$, define
\begin{gather*}
\quad
v_{k,(ij)}=\text{sign}\lb x_{ki}-x_{kj}\rb,\quad \th_{(ij)}=\frac{1}{\sqrt{M}}\lb v_{1,(ij)},v_{2,(ij)},\cdots,v_{p,{(ij)}}\rb^{\T},\\
\Thta=\lb
\mathbf{\theta}_{(12)},\cdots,\mathbf{\theta}_{(1n)},\mathbf{\theta}_{(23)},\cdots,\mathbf{\theta}_{(2n)},\cdots,\th_{(n-1,n)}\rb~,
\end{gather*}
where $M:=M(n)=n(n-1)/2$.
We can represent $\K_n$ in \eqref{eq: taukl}  as
\begin{equation}\label{eq:Kn}
\K_n=\sum_{1\leq i<j\leq n}\mathbf{\theta}_{(ij)}\mathbf{\theta}_{(ij)}^{\T}=\Thta\Thta^{\T}~.
\end{equation}
Notice that the $k$-th row of $\Thta$ contains information only related to
the $k$-th component of the original data ($k$-th row of $\X_n$),
thus rows of $\Thta$ are independent. We denote the $k$-th row of $\Thta$ by
\begin{align}\label{vk}
\vv_k^{\T}=\frac{1}{\sqrt M}\lb v_{k, (12)} , \cdots , v_{k, (1n)} ,  v_{k, (23)} ,  \cdots , v_{k, (2n)} , \cdots , v_{k, (n-1, n)}   \rb.
\end{align}
Further, if we look into the components of $\vv_k^{\T}$ carefully, we will
find out that not all of them are independent, e.g. $v_{k, (ij)}$ and
$v_{k, (i \ell)}$ are correlated when $j \neq \ell$.  To deal with such
dependence structure within $\vv_k^{\T}$,  a variation of Hoeffding decomposition is first introduced in \cite{LSD17} and further refined by \cite{BZG18}.  Specifically, let
\begin{align*}
v_{k, (i \cdot)}=\E \lb\text{sign}~(x_{ki}-x_{kj})~|~x_{ki}\rb,\quad v_{k, (\cdot j)}=\E \lb\text{sign}~(x_{ki}-x_{kj})~|~x_{kj}\rb~,
\end{align*}
 then  $v_{k, (i j)}$ can be decomposed into three parts
\begin{align*}
v_{k, (i j)}:= v_{k, (i \cdot)}+ v_{k, (\cdot j)}+\bar{v}_{k, (i j)}~,
\end{align*}
where $\{v_{k, (i \cdot)}, ~k=1, \cdots, p;~ i=1, \cdots, n\}$ are i.i.d.\ uniformly distributed on $[-1, 1]$ and
the three terms  $v_{k, (i \cdot)}, v_{k, (\cdot j)}$ and $\bar{v}_{k, (i j)}$ are pairwisely uncorrelated.
If we further set $$u_{k, (ij)}:=v_{k, (i \cdot)}+ v_{k, (\cdot j)}$$ and define
\begin{align}\label{ukvk}
\uu_k^{\T}&=\frac{1}{\sqrt M}\lb u_{k, (12)} , \cdots , u_{k, (1n)} ,  u_{k, (23)} ,  \cdots , u_{k, (2n)} , \cdots , u_{k, (n-1 n)}   \rb,\nonumber\\
\bar{\vv}_k^{\T}&=\frac{1}{\sqrt M}\lb \bar{v}_{k, (12)} , \cdots , \bar{v}_{k, (1n)} ,  \bar{v}_{k, (23)} ,  \cdots , \bar{v}_{k, (2n)} , \cdots , \bar{v}_{k, (n-1 n)}   \rb,
\end{align}
then the $k$-th row  of $\Thta$ can be expressed as the summation of two terms
\[
\vv_k^{\T}=\uu_k^{\T}+\bar{\vv}_k^{\T},
\]
with the following covariance structure
\[
\E \uu_k\uu^{\T}_k=\frac{1}{3M}\TT^\T\TT, \quad \E
\bar{\vv}_k\bar{\vv}^\T_k=\frac{1}{3M}\I_M~\quad \text{and}\quad \E
\vv_k\vv_k^\T=\frac{1}{3M}\Big(\TT^\T\TT+\I_M\Big)~.
\]
Here, $\TT$ is a $n \times M$ matrix with entries
\[
\TT=\big(t_{\ell, (ij)}\big)_{\ell, i<j},\quad t_{\ell, (ij)}:=\delta_{\ell
i}-\delta_{\ell j},\quad 1 \leq \ell\leq n, 1\leq i<j \leq n~.
\]
and $\lb\TT\TT^{\T}\rb^2=n(\TT\TT^{\T}),~\lb\TT^{\T}\TT\rb^2=n (\TT^{\T}\TT)$ (see more details in \cite{BZG18}).

\subsection{The main theorem}

Denote the eigenvalues of $\K_n$ as $\lambda_1,\cdots,\lambda_p$ in the
descending order. Of interest is the asymptotic behavior of
$\sum_{i=1}^pf(\lambda_j)/p$, the LSS of $\K_n$, where $f(x)$ is an analytic function on $[0, \infty)$. 
Note that
\begin{align}\label{lss}
\frac{1}{p}\sum_{i=1}^pf(\lambda_j) = \int f(x) \d F^{\K_n}(x).
\end{align}
As has already been established that the  ESD  of $\K_n$  converges in
probability to $F_c$, then asymptotically, the quantity \eqref{lss} will tend
to $\int f(x) \d F_c(x)$ almost surely. In order to study its second order
behavior, or the limiting distribution of  the normalized version of
$\sum_{i=1}^pf(\lambda_j)/p$, we define $F^{c_n}$ as the c.d.f, which is the
analogue for $F_c$ simply by replacing all the limiting value $c$  in the
density function given in \eqref{kendensity} with its finite sample counterpart
$c_n$. Its corresponding Stieltjes transform is thus denoted  as
$m_{F^{c_n}}(z)$. We will prove that the  convergence rate of $$\int f(x) \d
F^{\K_n}(x)-\int f(x) \d F^{c_n}(x)$$  is essentially $1/p$. To
this end, define
\[ G_n(x)=p\lb F^{\K_n}(x)- F^{c_n}(x)\rb.\]
Our main result is stated in the following theorem.

\begin{theorem}\label{thm:mainclt}
Suppose  $(p, n)\to \infty$,  $\x_1, \x_2,\cdots,\x_n$ is a random
sample from a $p$-dimensional population $\x$, where $\x$  has $p$ independent
components, all of which are absolutely continuous with respect to the Lebesgue
measure. Let $f_1,\cdots,f_k$ be functions on $\mathbb{R}$ and analytic on an
open interval containing the support of $d F_c(x)$ (defined in \eqref{kendensity}), then the random vector
\begin{equation*}
\lb \int f_1(x)\d G_n(x),\cdots, \int f_k(x)\d G_n(x)\rb
\end{equation*}
forms a tight sequence in $n$ and converges weakly to a Gaussian random vector $\lb X_{f_1}, \cdots, X_{f_k}\rb$, with means
\begin{align}\label{eq:EXf}
\E X_f=&-\frac{1}{2\pi i}\oint_\gamma f(z) \lb\frac{36c m^3(z)\lb
1+\frac{2}{3}c m(z) \rb}{\left[-9\lb 1+\frac{2}{3}c m(z)
\rb^2+4cm^2(z)\right]^2}\right.\\  \nonumber
&\quad-\left.\frac{2c^2m^3(z)\left[\lb 1+\frac{2}{3}cm(z)\rb^2+6+\frac{4}{3}cm(z)\right]}{-9\lb 1+\frac{2}{3}c m(z) \rb^2+4cm^2(z)}\rb\d z\\
&-\frac{1}{2\pi i}\oint_\gamma f(z) \frac{8c m^3(z)}{\lb 1+\frac{2}{3}c m(z)
\rb\left[-9\lb 1+\frac{2}{3}c m(z) \rb^2+4cm^2(z) \right]}\d z \nonumber
\end{align}
and covariance functions
\begin{align}\label{eq:VarXf}
&~\cov\lb X_{f_i}, X_{f_j}\rb =-\frac{1}{2\pi^2}\oiint f_i(z_i)f_j(z_j)
\frac{m'(z_i)m'(z_j)}{\lb m(z_i)-m(z_j)\rb^2}\d z_i\d z_j\\ \nonumber
& +\frac{1}{\pi^2}\oiint f_i(z_i)f_j(z_j)\frac{2cm'(z_i)m'(z_j)}{9\lb 1+\frac{2}{3}cm(z_i)\rb^2\lb 1+\frac{2}{3}cm(z_j)\rb^2}\d z_i \d z_j~,
\end{align}
where  $1\leq i\neq j\leq k$. The contours in \eqref{eq:EXf} and \eqref{eq:VarXf} (two in \eqref{eq:VarXf} are assumed to be non-overlapping) are closed and taken in the positive direction in the complex plane, each enclosing the support of $d F_c(x)$.
\end{theorem}

\begin{remark}
In Theorem~\ref{thm:mainclt}, independence  and continuity are the only two
assumptions imposed on the components of $\x$. There are no moment constraints  and furthermore, the $p$ components of $\x$ are not necessarily to be 
identically distributed.
\end{remark}

The general forms for the limiting means and covariance functions given in
\eqref{eq:EXf} and \eqref{eq:VarXf} are difficult to compute because they
involve both the variable $z$ and Stieltjes transform $m(z)$. The change of variable from $z$ to
$m(z)$ itself is very complicated,  not to mention that from the function
$f(z)$ to $m(z)$. For ease of computation, the following corollary provides an
alternative form of contour integrations for calculating the limiting means and
covariances, which only depends on a complex number $\xi$ that runs
counterclockwise along the unit circle.

\begin{corollary}\label{changee}
Under the same assumptions as in Theorem \ref{thm:mainclt}, the limiting means and covariance functions can be expressed as follows
\begin{align}\label{exf2}
\E X_f=&\lim_{r \downarrow 1}\frac{1}{2\pi i}\oint_{|\xi|=1} f\left(\frac 13+\frac 23 |1+\sqrt c \xi|^2\right)\cdot  \left(\frac{1}{(r^2\xi^2-1)\xi}-\frac{2}{\xi^3}+\frac{c}{2(\sqrt c+r \xi)^3}\right)d\xi\\
&+\lim_{r \downarrow 1}\frac{1}{2\pi i}\oint_{|\xi|=1} f\left(\frac 13+\frac 23
|1+\sqrt c \xi|^2\right)\cdot\left(\frac{3c}{(\sqrt c+r\xi)\xi^2}-\frac{c \sqrt
c}{(\sqrt c+r\xi)^2\xi^2} \right)d\xi\nonumber
\end{align}
and
\begin{align}\label{vex2}
\cov\lb X_{f_i}, X_{f_j}\rb =-\lim_{r \downarrow 1}\frac{1}{2\pi^2}\oiint  &f_i\left(\frac 13+\frac 23 |1+\sqrt c \xi_i|^2\right)f_j\left(\frac 13+\frac 23 |1+\sqrt c \xi_j|^2\right)\nonumber\\
& \cdot \left(\frac{1}{(\xi_i-r\xi_j)^2}-\frac{1}{\xi^2_i\xi^2_j}\right)d\xi_id\xi_j~.
\end{align}
\end{corollary}
The proof of this Corollary is relegated into  the supplementary file.

\begin{remark}
    Note that the diagonal elements of $\K_n$ are all one, which leads to the fact that $\tr(\K_n)=p$ is deterministic. Actually, by taking $f(x)=x$ and via some non-trivial calculations, we can show that, the centering term $\int x\d F^{c_n}(x)=p$, and the limiting mean $\E X_f$ and variance $\text{Var}(X_f)$ are both zero, which  further proves the validity of Theorem~\ref{thm:mainclt} in this degenerate case.
\end{remark}

\section{Application to Test of Independence}\label{appli}

Let $\x=\big( x_{1},x_{2},\cdots,x_{p}\big)^{\T}$ be a $p$-dimensional random
vector. Of interest is to  test
\begin{equation}\label{eq:H0}
 H_0: x_{1}, x_{2}, \cdots, x_{p} \text{~are~independent}
\end{equation}
based on $n$ samples under the  regime  $(p, n)\to \infty$.

\subsection{Test statistics and their limiting null distributions}\label{appl1}
Here we consider the following three test statistics based on the  Kendall's rank correlation matrix $\K_n$, 
\begin{gather*}
 Q_{\tau,2}=\tr(\K_n^2),\quad
Q_{\tau,4}=\tr(\K_n^4),\quad Q_{\tau,\log}=\log|\K_n|.
\end{gather*}
The first two test statistics $Q_{\tau,2}$ and $Q_{\tau,4}$ are polynomial functions (of order two and order four, respectively) of the  pairwise rank correlations among all the $p$ components and naturally we reject $H_0$ when their values are too large.
The third one  $Q_{\tau,\log}$ is motivated by the entropy loss  between the Kendall's rank correlation matrix $\K_n$ and its population version under $H_0$. To be more specific,  when all the components of $\x$ are independent, it is straightforward to verify that $\E \K_n=\I_p$. So the entropy loss between $\K_n$ and $\E \K_n$, which is defined as   
 $L(\K_n, \E \K_n)=\tr\K_n (\E \K_n)^{-1}-\log(|\K_n (\E \K_n)^{-1}|)-p$ \citep{james1961estimation, muirhead2009aspects, zhenghypothesis} reduces to $-\log| \K_n|$ under $H_0$, and we reject $H_0$ when the entropy loss is too large, or equivalently when $Q_{\tau,\log}$ is too small. 
 
On the other hand, all the three statistics $Q_{\tau,k}~(k=2,4)$ and $Q_{\tau,\log}$ can be directly linked to particular forms of LSS of $\K_n$ by taking  $f(x)=x^k,~k=2, 4$ and $f(x)=\log(x)$, respectively, i.e.
\begin{gather*}
Q_{\tau,k}=p\int x^k\d F^{\K_n}(x),~ k=2, 4\quad \text{and}\quad Q_{\tau,\log}=p\int\log(x)\d F^{\K_n}(x),
\end{gather*}
with their asymptotic fluctuation behaviors under $H_0$ fully
characterized by implementing the CLT for general LSS of $\K_n$ 
in Theorem \ref{thm:mainclt} or Corollary~\ref{changee}. Through some non-trivial calculations, 
the limiting  distributions for the three statistics under $H_0$ are given in following theorem.
\begin{theorem}\label{prop:indeptest}
	Assuming that conditions in Theorem~\ref{thm:mainclt} hold, under $H_0$, we have, as $(p,n)\rightarrow \infty$,
	\begin{gather*}
Q_{\tau,2}-p-\frac{4p^2}{9n}\xrightarrow{d} \mathcal{N} \lb \frac{14}{9}c^2-\frac{4}{9}c,~ \frac{64}{81}c^2\rb,\\
	Q_{\tau,4}-p-\frac{8p^2}{3n}-\frac{128p^3}{n^2}-\frac{16p^4}{81n^3}\xrightarrow{d} \mathcal{N} \lb \mu_{\tau,4},~ \sigma_{\tau,4}^2\rb,\\
		Q_{\tau,\log}+\frac{b}{a}\sqrt{pn}-(p+n)\log a +(n-p)\log(a-b\sqrt{p/n})\xrightarrow{d} \mathcal{N}(\mu_{\tau,\log},\sigma^2_{\tau,\log}),
		\end{gather*}
	where
	\begin{gather}\label{eq:sigma}
	\mu_{\tau,4}=-\frac{8}{3}c+\frac{140}{27}c^2+\frac{608}{81}c^3+\frac{112}{81}c^4,	\\
	\sigma_{\tau,4}^2=4\lb\frac{8}{3}c +\frac{352}{81}c^2+\frac{32}{27}c^3\rb^2+6\lb \frac{32}{27}c^{3/2}+\frac{64}{81}c^{5/2}\rb^2+\frac{2048}{81^2}c^4,\nonumber\\
	a=\frac{\sqrt{d_{c,+}}+\sqrt{d_{c,-}}}{2},~~b=\frac{\sqrt{d_{c,+}}-\sqrt{d_{c,-}}}{2},\nonumber\\
\mu_{\tau,\log}=-2\log a+\frac{1}{2}\log(a^2-b^2)+\log(a-b\sqrt{c})+\frac{2b}{a}\sqrt{c}-\frac{4ab\sqrt{c}-3cb^2}{4(a-b\sqrt{c})^2} \nonumber +\frac{b^2}{a^2}\\
+\left\{ \frac{3b^2-2ab\sqrt{c}}{4(a\sqrt{c}-b)^2}-\log\lb a-\frac{b}{\sqrt{c}}\rb\right\}\mathds{1}_{\{c>1\}}, \nonumber\\
\sigma^2_{\tau,\log}=2\log\frac{a^2}{a^2-b^2}-\frac{2b^2}{a^2}. \nonumber
	\end{gather}
\end{theorem}
The proof of this theorem is postponed to the supplementary file.

\begin{remark}
	\cite{Leung18} introduced three types of test statistics  that are constructed
	as sums or sums of squares of pairwise rank correlations,  including Kendall's
	$\tau$ as a special case. In fact, 
Theorem 4.1 in \cite{Leung18} shows that under $H_0$, when $p,n\rightarrow \infty$,
\begin{align} \label{nsh}
 \frac{n\left[\frac{1}{2}\tr\lb \K_n^2\rb-\frac{p}{2}-{p \choose 2}\mu_h\right]}{4p/9}\xrightarrow{d} \mathcal{N} (0, 1),
 \end{align}
where  $\mu_h=\frac{2(2n+5)}{9n(n-1)}=\frac{4}{9n}+O(n^{-2})$.
Under the high dimensional framework $(p, n)\to \infty$, \eqref{nsh} is consistent with our results for the limiting null distribution of $Q_{\tau,2}$ in Theorem~\ref{prop:indeptest}.
\end{remark}

  \subsection{Simulation experiments}
   In this section, we  conduct  numerical comparisons  to examine the finite sample
performance of the three proposed test statistics  $Q_{\tau,2}$, $Q_{\tau,4}$ and $Q_{\tau,\log}$ with some existing ones.
 Let $Z_{\alpha}$ be the upper-$\alpha$ quantile of the standard normal distribution
 at level $\alpha$. Based on Theorem~\ref{prop:indeptest}, we obtain three procedures
 for testing the null hypothesis in \eqref{eq:H0} as follows.
 \begin{gather*}
 \mbox{Reject } H_0  \mbox{ if } 
    \left\{Q_{\tau,2}-p-\frac{4p^2}{9n}> \frac{8p}{9n}Z_\alpha +\frac{14p^2}{9n^2}-\frac{4p}{9n}\right\}\\
 \mbox{ or }  \left\{Q_{\tau,4}-p-\frac{8p^2}{3n}-\frac{128p^3}{n^2}-\frac{16p^4}{81n^3}> \widehat{\sigma}_{\tau,4} Z_\alpha+\widehat{\mu}_{\tau,4}
 \right\} 
 \\
 \mbox{ or } \left\{ Q_{\tau,\log}+\frac{b}{a}\sqrt{pn}-(p+n)\log a +(n-p)\log(a-b\sqrt{p/n})< \widehat{\sigma}_{\tau,\log} Z_{1-\alpha}+\widehat{\mu}_{\tau,\log} 
  \right\},
 \end{gather*}
 where $\widehat{\mu}_{\tau,4},~\widehat{\sigma}_{\tau,4}^2,~\widehat{\mu}_{\tau,\log},~\widehat{\sigma}_{\tau,\log}^2$ are the ones by  replacing the limiting value $c$   in the terms $\mu_{\tau,4},~ \sigma_{\tau,4}^2,~\mu_{\tau,\log},~ \sigma_{\tau,\log}^2$ in \eqref{eq:sigma} with  its finite sample counterpart  $c_n=p/n$.

As for comparison, \citet{BZG18} proposed a test statistic $Q_{\tau,1}$, which is  based on the
largest eigenvalue $\lambda_1(\K_n)$ of $\K_n$.
They have shown that, under similar assumptions as in Theorem~\ref{thm:mainclt},
 \[
Q_{\tau,1}:= \frac{3}{2}n^{\frac{2}{3}}c_n^{\frac{1}{6}}d_{+,c_n}^{-\frac{2}{3}}\lb \lambda_1\lb \K_n\rb-\lambda_{+,c_n}\rb\xrightarrow{d} \text{TW}_1,
 \]
 where $d_{+,c_n}=(1+\sqrt{c_n})^2$, $\lambda_{+,c_n}=\frac{1}{3}+\frac{2}{3}d_{+,c_n}$ and $\text{TW}_1$ stands for the Tracy-Widom law of type I.
In addition to Kendall's rank correlation matrix model, there are some other
testing procedures  based on Spearman and Pearson-type correlation matrices,
denoted by $\S_n$ and $\R_n$, respectively. Here both $\S_n=\lb
s_{k\ell}\rb$ and $\R_n=\lb \rho_{k\ell}\rb$ are
$p\times p$ matrices where $s_{k\ell}$ and $\rho_{k\ell}$ are the Spearman and
Pearson correlation of the $k$-th and $\ell$-th row of $\X_n$ with
\begin{gather*}
s_{k\ell}=\frac{\sum_{i=1}^n\lb r_{ki}-\overline{r}_k \rb\lb r_{\ell i}-\overline{r}_\ell \rb}{\sqrt{\sum_{i=1}^n\lb r_{ki}-\overline{r}_k\rb^2}\sqrt{\sum_{i=1}^n\lb r_{\ell i}-\overline{r}_\ell \rb^2}},\quad  \overline{r}_k=\frac{1}{n}\sum_{i=1}^n r_{ki}=\frac{n+1}{2},\\
\rho_{k\ell}=\frac{\sum_{i=1}^n\lb x_{ki}-\overline{x}_k \rb\lb x_{\ell i}-\overline{x}_\ell \rb}{\sqrt{\sum_{i=1}^n\lb x_{ki}-\overline{x}_k\rb^2}\sqrt{\sum_{i=1}^n\lb x_{\ell i}-\overline{x}_\ell \rb^2}},\quad \overline{x}_k=\frac{1}{n}\sum_{i=1}^n x_{ki},
\end{gather*}
where $r_{ki}$
is the rank of $x_{ki}$ among $\lb x_{k1},\cdots, x_{kn}\rb$. Test statistics based on $\S_n$ and $\R_n$ include
\begin{itemize}
    \item [1.]
    $Q_{R,1}=\cfrac{n\lambda_1(\R_n)-( p^{1/2}+n^{1/2})^2}{(
    p^{1/2}+n^{1/2})(p^{-1/2}+n^{-1/2})^{1/3}}$,
    \citep{bao2012tracy,pillai2012edge};
    \item [2.]  $Q_{R,2}=\tr\lb\R_n\R_n^{\T}\rb-p-\frac{p^2}{n}$
    \citep{gao2017high};
    \item [3.]  $Q_{R,\max}=n \lb \max_{1\leq i<j\leq p}|\R_{ij}|\rb^2-4\log n+\log\log
    n$, \citep{jiang04};
    \item[4.]  $Q_{S,1}=n^{\frac{2}{3}}c_n^{\frac{1}{6}}d_{+,c_n}^{-\frac{2}{3}}\lb\lambda_1(\S_n) -d_{+,c_n}\rb$,
    \citep{bao19};
    \item [5.] $Q_{S,2}=\frac{n^2}{p^2}\tr\lb \S_n\S_n^{\T}\rb-\frac{n^2}{n-1}-\frac{n^2}{p}+\frac{n}{p}$,
    \citep{B15};
        \item [6.] $Q_{S,4}=\frac{n^4}{p^4}\tr\lb \S_n^4\rb-\frac{n^4}{(n-1)^3}-\frac{n^4}{p^3}-\frac{6n^4}{(n-1)p^2}-\frac{6n^4}{p(n-1)^2}$,
        \citep{B15};
    \item [7.]  $Q_{S,\max}=n \lb \max_{1\leq i<j\leq p}|\S_{ij}|\rb^2-4\log p+\log\log
    p$, \citep{zhou07}.
\end{itemize}

 To evaluate the finite sample performance of these test statistics, data are generated from various model scenarios  for different  $(p,n)$ combinations.
To examine Type I error rate, three models in the following are used with different  distributions for $\X_n=(x_{ij})_{p\times n}$.
\begin{itemize}
	\item[(I)]\textbf{ Mixed case:} $x_{ij}\sim \text{Gamma}(\text{shape}=4,\text{scale}=0.5)$ i.i.d. for $1\leq i\leq [p/2],~ 1\leq j\leq n$, $x_{ij}\sim t(5)$ i.i.d. for $[p/2]< i \leq p$, $1\leq j\leq n$;
	\item[(II)] \textbf{Mixed case:} $x_{ij}\sim \text{Cauchy}(\text{location}=0,\text{scale}=1)$ i.i.d. for $1\leq i\leq [p/2],~ 1\leq j\leq n$, $x_{ij}\sim t(5)$ i.i.d. for $[p/2]< i \leq p$, $1\leq j\leq n$;
	\item[(III)] \textbf{Heavy-tail case:} $x_{ij}\sim \text{Cauchy}(\text{location}=0,\text{scale}=1)$ i.i.d. for $1\leq i\leq p,~ 1\leq j\leq n$.
\end{itemize}
To examine their empirical power, both linear and non-linear alternatives are considered. A matrix $\Z_n$ with independent  components
is generated firstly following (I)$\sim$(III), then the data matrix $\X_n$ is constructed as follows.
\begin{itemize}
	\item[(IV)] \textbf{Toeplitz:}  $\X_n=\mathbf{A}\Z_n$, $\mathbf{A}=(a_{ij})_{p\times p}$, $a_{ii}=1$, for some $1<k_0<p$, $\forall 1\leq i\leq p$, $1\leq k\leq k_0$, $a_{i,i\pm k}=\rho_s^k$, $0<\rho_s<1$; $\forall 1\leq i\leq p$, $k_0< k<p$, $a_{i,i\pm k}=0$;
	\item[(V)] \textbf{Nonlinear correlation:}
	$x_{ij}=r_1z_{ij}+r_2z_{i+1,j}^2+r_3z_{i+2,j}^2+r_4e_{ij}$, $1\leq i\leq p,~1\leq j\leq n$ where $e_{ij}\sim N(0,1)$ i.i.d.
\end{itemize}
Here we set $k_0=[p/100]$, while for each distribution,  $\rho_s$ and $r_1\sim r_4$
are set differently to accommodate different degrees of dependence. All
empirical statistics are obtained using 1000 independent replicates. Note that
for Model (II) and (III), the elements of $\x$ do not have finite fourth order
moments, which fails to meet the requirement for Pearson-type test statistics,
therefore we eliminate all the corresponding results.

{\scriptsize
	\begin{table}[t]
		\caption{Empirical sizes for all test statistics under different distribution Models (I) $\sim $ (III).\label{Tab:SizeCom}}
		\resizebox{0.98\textwidth}{!}{
			\def\arraystretch{1.2}
			\begin{tabular}{ccc|ccc|ccc|ccc|cc}
				\hline
				\hlinewd{1.5pt}
				$p$ & $n$ & $c$ &  $Q_{\tau,\log}$& $Q_{\tau,4}$& $Q_{\tau,2}$ & $Q_{S,4}$ & $Q_{S,2}$ & $Q_{R,2}$ & $Q_{\tau,1}$ & $Q_{S,1}$ & $Q_{R,1}$ & $Q_{S,\max}$ & $Q_{R,\max}$
				\tabularnewline
				\hline
				\hlinewd{1.5pt}
			100  & 200  & 0.5  & 0.058 & 0.042 & 0.056  & 0.038  & 0.054  & 0.036  & 0.020  & 0.014  & 0.016  & 0.027  & 0.066\tabularnewline
			200  & 400  & 0.5 & 0.053 & 0.051  & 0.042  & 0.050  & 0.042  & 0.048  & 0.023  & 0.018  & 0.022  & 0.032  & 0.062\tabularnewline
			300  & 600  & 0.5  & 0.065 & 0.070  & 0.048  & 0.067  & 0.050  & 0.053  & 0.030  & 0.023  & 0.027  & 0.030  & 0.090\tabularnewline
			\hline 
			100  & 100  & 1  & 0.071 & 0.054  & 0.057  & 0.053  & 0.053  & 0.054  & 0.026  & 0.017  & 0.013  & 0.026  & 0.052\tabularnewline
			300  & 300  & 1 & 0.051 & 0.047  & 0.039  & 0.048  & 0.036  & 0.046  & 0.027  & 0.016  & 0.028  & 0.025  & 0.076\tabularnewline
			500  & 500  & 1  & 0.047 & 0.056  & 0.051  & 0.058  & 0.052  & 0.059  & 0.027  & 0.019  & 0.035  & 0.046  & 0.122\tabularnewline
			\hline 
			200  & 100  & 2  & 0.077& 0.068  & 0.054  & 0.064 & 0.049  & 0.057  & 0.033  & 0.017  & 0.018  & 0.017  & 0.050 \tabularnewline
			400  & 200  & 2  &0.054 & 0.045 & 0.049  & 0.036 & 0.042  & 0.044  & 0.038  & 0.021  & 0.025  & 0.022  & 0.135\tabularnewline
			600  & 300  & 2  &0.048  & 0.042 & 0.038  & 0.038 & 0.039  & 0.045  & 0.037  & 0.023  & 0.030  & 0.026  & 0.157\tabularnewline \hline
				\multicolumn{14}{l}{\textbf{For Distribution Model (I)}}  \tabularnewline
				\hlinewd{1.1pt}
			100  & 200  & 0.5 & 0.049 & 0.051 & 0.057  & 0.052  & 0.059  & -  & 0.020  & 0.009  & -  & 0.035  & -\tabularnewline
			200  & 400  & 0.5 & 0.058 & 0.048 & 0.050  & 0.045  & 0.050  & -  & 0.019  & 0.016  & -  & 0.039  & -\tabularnewline
			300  & 600  & 0.5 & 0.047 & 0.046 & 0.047  & 0.043  & 0.046  & -  & 0.027  & 0.023  & -  & 0.042  & -\tabularnewline
			\hline 
			100  & 100  & 1  & 0.069 & 0.060 & 0.060  & 0.050 & 0.053  & -  & 0.020  & 0.016  & -  & 0.025  & -\tabularnewline
			300  & 300  & 1  & 0.064 & 0.053  & 0.049  & 0.048  & 0.051  & -  & 0.025  & 0.018  & -  & 0.031  & -\tabularnewline
			500  & 500  & 1  & 0.045 & 0.055 & 0.053  & 0.054  & 0.052  & -  & 0.032  & 0.028  & -  & 0.039  & -\tabularnewline
			\hline 
			200  & 100  & 2  & 0.081 & 0.057 & 0.056  & 0.052  & 0.052  & -  & 0.036  & 0.022  & -  & 0.011  & -\tabularnewline
			400  & 200  & 2  & 0.047 & 0.061  & 0.052  & 0.055  & 0.047  & -  & 0.045  & 0.028  & -  & 0.039  & -\tabularnewline
			600  & 300  & 2  & 0.043 & 0.046  & 0.047  & 0.042  & 0.040  & -  & 0.037  & 0.023  & -  & 0.030  & -\tabularnewline 
			\hline
				\multicolumn{14}{l}{\textbf{For Distribution Model (II)}}  \tabularnewline 
				\hlinewd{1.1pt}
			100  & 200  & 0.5  & 0.057 & 0.055  & 0.044  & 0.052 & 0.040  & -  & 0.022  & 0.017  & -  & 0.036  & -\tabularnewline
			200  & 400  & 0.5  & 0.047  & 0.048 & 0.048  & 0.048  & 0.047  & -  & 0.022  & 0.017  & -  & 0.039  & -\tabularnewline
			300  & 600  & 0.5  & 0.055 & 0.050 & 0.057  & 0.049  & 0.054  & -  & 0.023  & 0.017  & -  & 0.050  & -\tabularnewline
			\hline 
			100  & 100  & 1  & 0.065  & 0.054 & 0.055  & 0.048 & 0.050  & -  & 0.019  & 0.013  & -  & 0.018  & -\tabularnewline
			300  & 300  & 1  & 0.059 & 0.056  & 0.040  & 0.056  & 0.039  & -  & 0.035  & 0.024  & -  & 0.038  & -\tabularnewline
			500  & 500  & 1  & 0.062 & 0.051 & 0.056  & 0.048  & 0.053  & -  & 0.026  & 0.021  & -  & 0.039  & -\tabularnewline
			\hline 
			200  & 100  & 2  & 0.076 & 0.059  & 0.056  & 0.052  & 0.047  & -  & 0.052  & 0.037  & -  & 0.027  & -\tabularnewline
			400  & 200  & 2  &  0.052 & 0.053 & 0.050  & 0.053 & 0.048  & -  & 0.040  & 0.023  & -  & 0.020  & -\tabularnewline
			600  & 300  & 2  & 0.044 & 0.066 & 0.061  & 0.063  & 0.059  & -  & 0.044  & 0.033  & -  & 0.050  & -\tabularnewline
			\hline
				\multicolumn{14}{l}{\textbf{For Distribution Model (III)}}  \tabularnewline
				\hline
				\hlinewd{1.1pt}
			\end{tabular}}
		\end{table}
	}

{\scriptsize
	\begin{table}[t]
		\caption{Empirical power for all test statistics under Toeplitz population matrix Model (IV) with different distributions.\label{Tab:PowerCom5}}
		\resizebox{0.98\textwidth}{!}{
			\def\arraystretch{1.2}
			\begin{tabular}{ccc|ccc|ccc|ccc|cc}
				\hline
				\hlinewd{1.5pt}
				$p$ & $n$ & $c$ & $Q_{\tau,\log}$ & $Q_{\tau,4}$ &$Q_{\tau,2}$ &  $ Q_{S,4}$ & $Q_{S,2}$ & $Q_{R,2}$ & $Q_{\tau,1}$ & $Q_{S,1}$ & $Q_{R,1}$ & $Q_{S,\max}$ & $Q_{R,\max}$\tabularnewline
				\hlinewd{1.2pt}
				\hline
100  & 200  & 0.5  & 0.871  & 0.863  & 0.898  & 0.844  & 0.894  & 0.880  & 0.190  & 0.163  & 0.168  & 0.246  & 0.290\tabularnewline
200  & 400  & 0.5  & 1  & 1  & 1  & 1  & 1  & 1  & 0.494  & 0.452  & 0.405  & 0.885  & 0.850\tabularnewline
300  & 600  & 0.5  & 1  & 1  & 1  & 1  & 1  & 1  & 0.768  & 0.732  & 0.652  & 1  & 1\tabularnewline
\hline 
100  & 100  & 1  & 0.424  & 0.410  & 0.405  & 0.373  & 0.391  & 0.368  & 0.108  & 0.066  & 0.069  & 0.058  & 0.075\tabularnewline
300  & 300  & 1  & 0.991  & 0.995  & 0.994  & 0.995  & 0.992  & 0.996  & 0.332  & 0.283  & 0.265  & 0.488  & 0.511\tabularnewline
500  & 500  & 1  & 1  & 1  & 1  & 1  & 1  & 1  & 0.584  & 0.534  & 0.504  & 0.987  & 0.981\tabularnewline
\hline 
200  & 100  & 2  & 0.419 & 0.401  & 0.470  & 0.366  & 0.437  & 0.419  & 0.112  & 0.055  & 0.056  & 0.038  & 0.097\tabularnewline
400  & 200  & 2  & 0.865  & 0.915  & 0.908  & 0.892  & 0.909  & 0.870  & 0.245  & 0.182  & 0.155  & 0.178  & 0.237\tabularnewline
600  & 300  & 2  &0.991  & 0.993 & 0.996  & 0.991  & 0.995  & 0.994  & 0.369  & 0.287  & 0.241  & 0.437  & 0.484\tabularnewline
\hline
				\multicolumn{14}{l}{$\rho_{s}=0.06$ {\bf for Distribution Model (I)}}  \tabularnewline
				\hlinewd{1.1pt}
100  & 200  & 0.5  & 0.889  & 0.892  & 0.899  & 0.868  & 0.882  & -  & 0.225  & 0.187  & -  & 0.466  & -\tabularnewline
200  & 400  & 0.5  & 1  & 1  & 1  & 1  & 1  & -  & 0.583  & 0.528  & -  & 0.999  & -\tabularnewline
300  & 600  & 0.5  & 1  & 1  & 1  & 1  & 1  & -  & 0.838  & 0.804  & -  & 1  & -\tabularnewline
\hline 
100  & 100  & 1  & 0.464  & 0.447  & 0.473  & 0.405  & 0.443  & -  & 0.095  & 0.072  & -  & 0.086  & -\tabularnewline
300  & 300  & 1  & 1  & 0.998  & 0.999  & 0.997  & 0.996  & -  & 0.391  & 0.336  & -  & 0.936  & -\tabularnewline
500  & 500  & 1  & 1  & 1  & 1  & 1  & 1  & -  & 0.683  & 0.631  & -  & 1  & -\tabularnewline
\hline 
200  & 100  & 2  & 0.508 & 0.450  & 0.494  & 0.405  & 0.448  & -  & 0.137  & 0.083  & -  & 0.099  & -\tabularnewline
400  & 200  & 2  &0.901  & 0.922  & 0.925  & 0.900  & 0.912  & -  & 0.239  & 0.165  & -  & 0.457  & -\tabularnewline
600  & 300  & 2  & 0.995 & 0.999 & 0.999  & 0.998  & 0.999  & -  & 0.399  & 0.305  & -  & 0.924  & -\tabularnewline
\hline
				\multicolumn{14}{l}{$\rho_{s}=0.03$ {\bf for Distribution Model (II)}}  \tabularnewline
				\hlinewd{1.1pt}
100  & 200  & 0.5  & 0.843  & 0.818  & 0.833  & 0.786  & 0.806  & -  & 0.181  & 0.154  & -  & 0.233  & -\tabularnewline
200  & 400  & 0.5  & 1  & 1  & 1  & 1  & 1  & -  & 0.461  & 0.427  & -  & 0.823  & -\tabularnewline
300  & 600  & 0.5  & 1  & 1  & 1  & 1  & 1  & -  & 0.743  & 0.708  & -  & 0.999  & -\tabularnewline
\hline 
100  & 100  & 1  & 0.410  & 0.397  & 0.390  & 0.358  & 0.353  & -  & 0.071  & 0.049  & -  & 0.053  & -\tabularnewline
300  & 300  & 1  & 0.985  & 0.993  & 0.993  & 0.990  & 0.992  & -  & 0.337  & 0.274  & -  & 0.508  & -\tabularnewline
500  & 500  & 1  & 1  & 1  & 1  & 1  & 1  & -  & 0.592  & 0.524  & -  & 0.978  & -\tabularnewline
\hline 
200  & 100  & 2  &  0.501 & 0.460  & 0.474  & 0.411  & 0.425  & -  & 0.110  & 0.075  & -  & 0.034  & -\tabularnewline
400  & 200  & 2  & 0.873  & 0.874  & 0.899  & 0.848  & 0.882  & -  & 0.204  & 0.154  & -  & 0.146  & -\tabularnewline
600  & 300  & 2  & 0.984 & 0.990 & 0.994  & 0.988  & 0.992  & -  & 0.340  & 0.277  & -  & 0.421  & -\tabularnewline
\hline
				\multicolumn{14}{l}{{\bf $\rho_{s}=0.02$ for Distribution Model (III) }}  \tabularnewline
				\hlinewd{1.1pt}
			\end{tabular}
		}
	\end{table}}

		{\scriptsize
		\begin{table}[t]
			\caption{Empirical power for all feasible test statistics under nonlinear correlation Model (V) with heavy tailed distributions (II) and (III).\label{Tab:PowerCom7}}
			\resizebox{0.9\textwidth}{!}{
				\def\arraystretch{1.2}
				\begin{tabular}{ccc|ccc|cc|ccc}
					\hline
					\hlinewd{1.5pt}
					$p$ & $n$ & $c$ &  $Q_{\tau,\log}$& $Q_{\tau,4}$ & $Q_{\tau,2}$ & $Q_{S,4}$ & $Q_{S,2}$ & $Q_{\tau,1}$ & $Q_{S,1}$  & $Q_{S,\max}$ \tabularnewline
					\hline
					\hlinewd{1.3pt}
100 & 200 & 0.5 & 0.949 & 0.931 & 0.947 & 0.899 & 0.929 & 0.235 & 0.179 & 0.765\tabularnewline
200 & 400 & 0.5 & 1 & 1 & 1 & 1 & 1 & 0.551 & 0.491 & 1\tabularnewline
300 & 600 & 0.5 & 1 & 1 & 1 & 1 & 1 & 0.804 & 0.750 & 1\tabularnewline
\hline 
100 & 100 & 1 & 0.573 & 0.495 & 0.561 & 0.433 & 0.478 & 0.091 & 0.065 & 0.180\tabularnewline
300 & 300 & 1 & 1 & 0.994 & 0.999 & 0.992 & 0.996 & 0.374 & 0.313 & 0.997\tabularnewline
500 & 500 & 1 & 1 & 1 & 1 & 1 & 1 & 0.654 & 0.584 & 1\tabularnewline
\hline 
200 & 100 & 2 & 0.641 & 0.485 & 0.545 & 0.424 & 0.456 & 0.123 & 0.071 & 0.132\tabularnewline
400 & 200 & 2 & 0.962& 0.934 & 0.963 & 0.896 & 0.934 & 0.268 & 0.187 & 0.725\tabularnewline
600 & 300 & 2 &  0.998& 0.999 & 0.999 & 0.999 & 0.999 & 0.350 & 0.263 & 0.998\tabularnewline
\hline
					\multicolumn{10}{l}{$r_1=0.01,~r_2=0.02,~r_3=0.006,~r_4=0.5$ \textbf{for Distribution Model (II)}}  \tabularnewline \hline
					\hlinewd{1.1pt}
100 & 200 & 0.5 & 0.699 & 0.667 & 0.708 & 0.620 & 0.670 & 0.119 & 0.092 & 0.159\tabularnewline
200 & 400 & 0.5 & 0.994 & 0.992 & 0.995 & 0.987 & 0.992 & 0.260 & 0.219 & 0.580\tabularnewline
300 & 600 & 0.5 & 1 & 1 & 1 & 1 & 1 & 0.482 & 0.437 & 0.953\tabularnewline
\hline 
100 & 100 & 1 & 0.383 & 0.318 & 0.349 & 0.271 & 0.310 & 0.066 & 0.034 & 0.041\tabularnewline
300 & 300 & 1 & 0.935 & 0.926 & 0.946 & 0.910 & 0.932 & 0.192 & 0.146 & 0.279\tabularnewline
500 & 500 & 1 & 1 & 1 & 1 & 1 & 1 & 0.335 & 0.286 & 0.801\tabularnewline
\hline 
200 & 100 & 2 & 0.407 & 0.310 & 0.342 & 0.266 & 0.284 & 0.095 & 0.062 & 0.031\tabularnewline
400 & 200 & 2 &0.730  & 0.665 & 0.717 & 0.611 & 0.656 & 0.130 & 0.083 & 0.097\tabularnewline
600 & 300 & 2 &0.944  & 0.930 & 0.954 & 0.908 & 0.940 & 0.210 & 0.161 & 0.219\tabularnewline \hline
					\multicolumn{10}{l}{$r_1=0.002,~r_2=0.005,~r_3=0.0015,~r_4=0.5$ \textbf{for Distribution Model (III)}} \tabularnewline
						\hlinewd{1.1pt}
			\end{tabular}}
	\end{table}}

Table~\ref{Tab:SizeCom} shows empirical sizes of all the test
statistics  under different distribution models (I)$\sim$(III).   It can be seen
that for the nominal level $\alpha=5\%$, all the Frobenius-norm-type test
statistics ($Q_{\tau,2}$, $Q_{\tau,4}$, $Q_{S,4}$, $Q_{S,2}$, $Q_{R,2}$) have very accurate sizes close to $5\%$, while spectral-norm-type ($Q_{\tau,1}$, $Q_{S,1}$, $Q_{S,1}$) and maximum-norm-type ($Q_{S,\max}$, $Q_{R,\max}$) test statistics are a little bit undersized.  Although
such bias shrinks when the sample size becomes larger, it still exists even for
very large $n$ due to the slow convergence of extreme eigenvalues to the
Tracy-Widom distribution. Moreover, $Q_{R,1}$, $Q_{R,2}$ and $Q_{R,\max}$ work only for  model (I) due to moment restrictions. It's observed that our test statistic $Q_{\tau,\log}$ is a little oversized in the $p>n$ cases when $p,n$ are relatively small. However such bias significantly reduces when $p,n$ increase. 

As for the linear correlated alternatives,
Table~\ref{Tab:PowerCom5} presents the empirical powers of all test statistics
under the Toeplitz population matrix Model (IV) for various distributions. It can be seen that our test statistics $Q_{\tau,\log}$ and $Q_{\tau,2}$
 perform best under Model (IV). 
 Moreover, among all the three Frobenius-norm-type statistics, our test statistic $Q_{\tau,2}$ demonstrates superiority over the other two
across all scenarios and $(p,n)$ combinations. $Q_{S,2}$ always has inferior
power and $Q_{R,2}$ doesn't work for heavy-tailed distributions.
As for the nonlinear correlated structure for Model (V), empirical powers of all the test statistics are shown in Table~\ref{Tab:PowerCom7}. It's obvious that our proposed test statistics  $Q_{\tau,2}$, $Q_{\tau,4}$ and $Q_{\tau,\log}$ all demonstrate significant superiority over all the others  under heavy tailed distributions.

\section{Proofs of Theorem \ref{thm:mainclt}}\label{mainre}
Generally speaking, the proof of our main result follows similar routine as establishing the CLT for LSS of a large dimensional  sample covariance matrix given in \citet{BS04}. However, the model we considered here is much more complicated than the well studied sample covariance matrix model as has already been discussed in the introduction. Extra new techniques are thus needed to overcome such difficulties.
To this end, we have developed three  lemmas (Lemma \ref{lem:generalvkbound}, Lemma \ref{lem:vkbound} and Lemma \ref{lem:quadratics}), all considering the expectation  of the product of a sequence of quadratic forms. Actually,  such results  serve as the cornerstone for proving our main result.

\subsection{Sketch of the proof of Theorem \ref{thm:mainclt}}

Let $\gamma$ be a closed contour enclosing the support of $d F_c(x)$, taken in the positive direction in the complex plane, then from Cauchy's integral formula, we have
\begin{equation*}
\int f(x) \d G_n(x)=-\frac{1}{2\pi i}\oint_{\mathcal{C}} f(z)M_n(z)\d z~,
\end{equation*}
where
$M_n(z)=p\lb m_{F^{\K_n}}(z)-m_{F^{c_n}}(z)\rb$. The proof of our main result Theorem~\ref{thm:mainclt} is formulated as to show the convergence of random process $M_n(z)$. More precisely, $M_n(z)$ is a random two-dimensional process defined on the contour $\mathcal{C}$ of the complex plane. It can be viewed as a random element in the metric space $C\lb \mathcal{C}, \mathbb{R}^2\rb$ of continuous functions from $\mathcal{C}$  to $\mathbb{R}^2$. Our target is to prove the weak convergence of $M_n(z)$ and the result is stated in the following lemma.
\begin{lemma}\label{lem:main}
    Under the same assumptions as in Theorem \ref{thm:mainclt}, $M_n(z)$ forms a tight sequence on $\mathcal{C}$ and converges weakly to a two-dimensional Gaussian process $M(z)$ satisfying, for $z\in \mathcal{C}$,
    \begin{align}\label{eq:lemMean}\nonumber
    \E M(z)=&~\frac{36c m_{F_{c}}^3(z)\lb 1+\frac{2}{3}c m_{F_{c}}(z) \rb}{\left[-9\lb 1+\frac{2}{3}c m_{F_{c}}(z) \rb^2+4cm_{F_{c}}^2(z)\right]^2}-\frac{2c^2m_{F_{c}}^3(z)\left[\lb 1+\frac{2}{3}cm_{F_{c}}(z)\rb^2+6+\frac{4}{3}cm_{F_{c}}(z)\right]}{-9\lb 1+\frac{2}{3}c m_{F_{c}}(z) \rb^2+4cm_{F_{c}}^2(z)}\\
    &\quad+\frac{8c m_{F_{c}}^3(z)}{\lb 1+\frac{2}{3}c m_{F_{c}}(z) \rb\left[-9\lb 1+\frac{2}{3}c m_{F_{c}}(z) \rb^2+4cm_{F_{c}}^2(z) \right]}~
    \end{align}
    and for $z_1,z_2\in \mathcal{C}$,
    \begin{equation}\label{eq:lemVar}
    \cov\lb M(z_1), M(z_2)\rb=
    \frac{2m_{F_{c}}'(z_1)m_{F_{c}}'(z_2)}{\lb m_{F_{c}}(z_1)-m_{F_{c}}(z_2)\rb^2}-\frac{2}{(z_1-z_2)^2}-\frac{8cm'_{F_{c}}(z_1)m'_{F_{c}}(z_2)}{9\lb 1+\frac{2}{3}cm_{F_{c}}(z_1)\rb^2\lb 1+\frac{2}{3}cm_{F_{c}}(z_2)\rb^2}.\\
    \end{equation}
\end{lemma}


\subsection{Proof of Lemma \ref{lem:main}}\label{pse}
Firstly we decompose $M_n(z)$ into the summation of a random  part $M_n^1(z)$ and a determinist part $M_n^2(z)$,  where
\begin{equation*}
M_n^1(z)=p\lb m_{F^{\K_n}}(z)-\E m_{F^{\K_n}}(z) \rb\quad\text{and}\quad M^2_n(z)=p\lb\E m_{F^{\K_n}}(z)-m_{F^{c_n}}(z)\rb.
\end{equation*}
The proof of Lemma \ref{lem:main} is then complete if we can verify the following three steps:
\begin{itemize}[leftmargin=18mm, itemsep=-10pt]
\item [Step 1:] Finite dimensional convergence of $M_n^1(z)$ in distribution to a centered multivariate Gaussian random vector with covariance function given by \eqref{eq:lemVar};\\
\item [Step 2:] Tightness of  $M_n^1(z)$;\\
\item [Step 3:]  Convergence of $M_n^2(z)$ to the mean function given by \eqref{eq:lemMean}.
\end{itemize}

\subsubsection{Step 1: Finite dimensional convergence of $M_n^1(z)$ in distribution}

In the first step, we will show that, for any complex numbers $z_1 \cdots z_r \in \mathcal{C}^+$, the $r$ dimensional random vector $(M^1_n(z_1)\cdots M^1_n(z_r))^{\T}$ is jointly Gaussian. To this end, we will show that
the sum
$$\sum_{i=1}^r\alpha_iM_n^1(z_i),\quad z_i\in \mathcal{C}^+$$
forms a tight sequence of random functions for $z_i$ and will converge in distribution to a Gaussian random variable.
The main strategy of the proof is based on the martingale CLT  given in Lemma~\ref{lem:bili08}, together with our newly established Lemma \ref{lem:generalvkbound}$\sim$\ref{lem:quadratics}.

Define
\begin{align*}
\mathbf{\Thta}_k&=\lb\vv_1,\cdots,\vv_{k-1},\vv_{k+1},\cdots,\vv_p \rb^{\T},
\quad \K_{n,k}=\mathbf{\Thta}_k\mathbf{\Thta}_k^{\T}~,\\
\A_k&=\mathbf{\Thta}_k^{\T}\lb\K_{n,k}-z\I_{p-1}\rb^{-2}\mathbf{\Thta}_k~,\quad
\B_k=\mathbf{\Thta}_k^{\T}\lb\K_{n,k}-z\I_{p-1}\rb^{-1}\mathbf{\Thta}_k~,\\
\beta_k&=\cfrac{1}{-\vv_k^{\T}\vv_k+z+\vv_k^{\T}\B_k\vv_k}~,\\
\widetilde{b}_k&=\cfrac{1}{-1+z+\frac{1}{3M}\tr\lb \TT\B_k\TT^{\T}\rb+\frac{1}{3M}\tr\lb \B_k\rb}~,\\
b_n&=\cfrac{1}{-1+z+\frac{1}{3M}\E\tr\lb \TT\B_k\TT^{\T}\rb+\frac{1}{3M}\E\tr(\B_k)}~,\\
g_k&=\vv_k^{\T}\vv_k-1-\vv_k^{\T}\B_k\vv_k+\frac{1}{3M}\tr\lb \TT\B_k\TT^{\T}\rb+\frac{1}{3M}\tr \B_k~,\\
h_k&=\vv_k^{\T}\A_k\vv_k-\frac{1}{3M}\tr\lb \TT\A_k\TT^{\T}\rb-\frac{1}{3M}\tr \A_k~,\\
a_k&=-g_k\beta_k\widetilde{b}_k\lb 1+\vv_k^{\T}\A_k\vv_k\rb~,\\
d_k&=h_k\widetilde{b}_k~.
\end{align*}
Let $\E_k(\cdot)$ denote  the conditional expectation with respect to the $\sigma-$field generated by $\vv_1,\cdots,\vv_k$, we have
\begin{align}\label{mn1}
M_n^1(z)&=p\lb m_{F^{\K_n}}(z)-\E m_{F^{\K_n}}(z) \rb\nonumber\\
&=\sum_{k=1}^p\lb\E_k-\E_{k-1}\rb\left\{\tr\lb\K_n-z\I_p\rb^{-1}-\tr\lb\K_{n,k}-z\I_{p-1}\rb^{-1}\right\}\nonumber\\
&=\sum_{k=1}^p\lb\E_k-\E_{k-1}\rb\cfrac{1+\vv_k^{\T}\mathbf{\Thta}_k^{\T}\lb \K_{n,k}-z\I_{p-1}\rb^{-2}\mathbf{\Thta}_k\vv_k}{\vv_k^{\T}\vv_k-z-\vv_k^{\T}\mathbf{\Thta}_k^{\T}\lb \K_{n,k}-z\I_{p-1}\rb^{-1}\mathbf{\Thta}_k\vv_k}\nonumber\\
&=\sum_{k=1}^p\lb \E_k-\E_{k-1}\rb \left\{a_k-d_k+\cfrac{1+\frac{1}{3M}\tr\lb \TT\A_k\TT^{\T}\rb+\frac{1}{3M}\tr\A_k}{1-z-\frac{1}{3M}\tr\lb \TT \B_k\TT^{\T}\rb-\frac{1}{3M}\tr\B_k}\right\}\nonumber\\
&=\sum_{k=1}^p\lb \E_k-\E_{k-1}\rb (a_k-d_k)=\sum_{k=1}^p\lb \E_k-\E_{k-1}\rb a_k- \E_k d_k~.
\end{align}
By applying the equality $$\beta_k=\widetilde{b}_k+\beta_k\widetilde{b}_kg_k~,$$ we have
\begin{align*}
a_k = &-g_k\beta_k\widetilde{b}_k\lb 1+\vv_k^{\T}\A_k\vv_k\rb~\\
=&~-\widetilde{b}_k^2g_k\lb 1+\frac{1}{3M}\tr\lb \TT\A_k\TT^{\T}\rb+\frac{1}{3M}\tr \A_k\rb-h_kg_k\widetilde{b}_k^2-\beta_k\widetilde{b}_k^2\lb 1+\vv_k^{\T}\A_k\vv_k\rb g_k^2\\
:=&~a_{k1}+a_{k2}+a_{k3},
\end{align*}
which together with \eqref{mn1} implies that
\[M_n^1(z)=\sum_{k=1}^p\lb\E_k-\E_{k-1}\rb\lb a_{k1}+a_{k2}+a_{k3}\rb-\E_kd_k.\]
Next, we will show that the contribution of $a_{k2}$ and $a_{k3}$ to $M_n^1(z)$ can be negligible as $n\to \infty$. For the contribution of $a_{k3}$, we have
\begin{align*}
\E\left\lvert\sum_{k=1}^p \lb\E_k-\E_{k-1}\rb a_{k3}\right\rvert^2 &\leq C_0\sum_{k=1}^p\E\lb 1+\vv_k^{\T}\A_k\vv_k\rb^2g_k^4\\
&\leq C_0 \sum_{k=1}^p\left\{\E\Big[ \Big( 1+\frac{1}{3M}\tr\lb \TT\A_k\TT^{\T}\rb+\frac{1}{3M}\tr \A_k\Big)^2g_k^4\Big]+\E \lb h_k^2g_k^4\rb \right\}.
\end{align*}
Using the Hoeffding decomposition  for $\vv_k$ which is   denoted as  $\vv_k=\uu_k+\vb_k$ (see \eqref{ukvk})  together with  Proposition 2.2 and 2.4 in \cite{BZG18},  we have
\begin{align}\label{gk4}
~\E g_k^4&=~\E \left[\vv_k^{\T}\lb \I_M-\B_k\rb\vv_k+\frac{1}{3M}\tr\lb \TT\B_k\TT^{\T}\rb+\frac{1}{3M}\tr\B_k-1\right]^4\nonumber\\
&=~\E\left[ \uu_k^{\T}\lb \I_M-\B_k\rb\uu_k-\frac{1}{3M}\tr\lb \TT\lb \I_M-\B_k\rb\TT^{\T}\rb+2\vb_k^{\T}\lb \I_M-\B_k\rb\uu_k\right.\nonumber\\
&\left.\quad~+~\vb_k^{\T}\lb \I_M-\B_k\rb\vb_k-\frac{1}{3M}\tr\lb\I_M-\B_k\rb\right]^4\nonumber\\
&\leq C_0\left\{\E\lb \uu_k^{\T}\lb \I_M-\B_k\rb\uu_k-\frac{1}{3M}\tr\lb \TT\lb \I_M-\B_k\rb\TT^{\T}\rb \rb^4+\E\lb \vb_k^{\T}\lb \I_M-\B_k\rb\uu_k\rb^4\right.\nonumber\\
&\left.\quad~+~\E\lb \vb_k^{\T}\lb \I_M-\B_k\rb\vb_k-\frac{1}{3M}\tr\lb \I_M-\B_k\rb \rb^4\right\}\nonumber\\
&\prec\lb\frac{\tr\lv \TT\lb \I_M-\B_k\rb\TT^{\T}\rv^2}{9M^2}\rb^2+\lb \frac{n}{M^2}\tr \lv \I_M-\B_k\rv^2\rb^2+\lb \frac{1}{M^2}\sum_{l=1}^n\Big|\sum_{j=l+1}^n\lb\TT\lb \I_M-\B_k\rb\rb_{j,(lj)}\Big|^2\rb^2\nonumber\\
&=O(n^{-2}),
\end{align}
which implies that
\begin{align}\label{eegk}
\E\left[ \lb 1+\frac{1}{3M}\tr\lb \TT\A_k\TT^{\T}\rb+\frac{1}{3M}\tr \A_k\rb^2g_k^4\right]=O(n^{-2}).
\end{align}

According to Lemma~\ref{lem:vkbound},
\begin{align}\label{egk}
\E h_k^2g_k^4\leq \sqrt{\E h_k^4\cdot \E g_k^8}\leq \sqrt{O(n^{-2})\cdot O(n^{-4})}=O(n^{-3}).
\end{align}
Collecting \eqref{eegk} and \eqref{egk}, we have
\[\E\Big\rvert\sum_{k=1}^p \lb\E_k-\E_{k-1}\rb a_{k3}\Big\rvert^2\leq C_0\sum_{k=1}^p O(\frac{1}{n^2})=O(n^{-1}).\]
Similarly,
\begin{align*}
&~\E\Big\lvert\sum_{k=1}^p \lb\E_k-\E_{k-1}\rb a_{k2}\Big\rvert^2=\sum_{k=1}^{p}\E\left\lvert\lb\E_k-\E_{k-1}\rb a_{k2}\right\rvert^2\\
\leq &~\sum_{k=1}^{p}\E \Big\lvert h_kg_k\widetilde{b}_k^2\Big\rvert^2\leq C_0\sum_{k=1}^{p}\E~ h_k^2g_k^2\leq  C_0\sum_{k=1}^{p}\sqrt{\E h_k^4\cdot\E g_k^4}\\
=&~C_0\sum_{k=1}^p \sqrt{O(n^{-2})\cdot O(n^{-2})}=O(n^{-1}),
\end{align*}
which implies that
\begin{align*}
M_n^1(z)
&=\sum_{k=1}^p\lb\E_k-\E_{k-1}\rb\lb a_{k1}+a_{k2}+a_{k3}\rb-\E_kd_k
=\sum_{k=1}^{p}\E_k\lb a_{k1}-d_k\rb+o(1)\\
&=\sum_{k=1}^p \E_k \left\{-\widetilde{b}_k^2g_k\Big( 1+\frac{1}{3M}\tr( \TT\A_k\TT^{\T})+\frac{1}{3M}\tr\A_k\Big)-h_k\widetilde{b}_k\right\}+o(1)\\
&=\sum_{k=1}^p \E_k\left(\frac{\d~}{\d z}\widetilde{b}_k(z)g_k(z)\right)+o(1)~,
\end{align*}
where $\E_k\frac{\d{~}}{\d z}\lb \widetilde{b}_k(z)g_k(z)\rb$ is a sequence of martingale difference.
Since
\[
\sum_{i=1}^r\alpha_iM_n^1(z_i)=\sum_{k=1}^{p}\sum_{i=1}^r\alpha_i\E_k\frac{\d{~}}{\d z_i}\lb \widetilde{b}_k(z_i)g_k(z_i)\rb+o(1):=\sum_{k=1}^{p}\lb\sum_{i=1}^r\alpha_i Y_k(z_i)\rb+o(1)~,
\]
by applying the martingale central limit theorem (\citet{bili08})  as stated in Lemma~\ref{lem:bili08}, it is enough to verify
\begin{equation}\label{eq:MCLTC1}
\sum_{k=1}^p\E\lb \Big(\sum_{i=1}^r\alpha_i Y_k(z_i)\Big)^2\mathds{1}_{\left\{\left\lvert\sum_{i=1}^r\alpha_i Y_k(z_i)\right\rvert\geq \veps\right\}}\rb\xrightarrow{p} 0,
\end{equation}
and prove that under the assumptions of our Theorem \ref{thm:mainclt}, for $z_1, z_2 \in \mathbb{C}^{+}$, the term
\begin{equation}\label{limito}
\sum_{k=1}^{p}\E_{k-1}\Big( Y_k(z_1) Y_k(z_2)\Big)
\end{equation}
converges in probability to a constant, or to determine the limit of
\begin{equation}\label{eq:MCLTC2}
\sum_{k=1}^{p}\E_{k-1}\left[\E_k\lb \widetilde{b}_k(z_1)g_k(z_1)\rb \E_k\lb \widetilde{b}_k(z_2)g_k(z_2)\rb\right].
\end{equation}
Actually, we have
$$\frac{d}{dz_1}\frac{d}{dz_2}\eqref{eq:MCLTC2}=\eqref{limito}~.$$
As for the condition \eqref{eq:MCLTC1}, by taking \eqref{gk4} and \eqref{egk} into account, we have
\begin{align*}
&~\E\left\lvert ~\E_{k}\frac{\d~}{\d z}\lb \widetilde{b}_k(z)g_k(z)\rb~\right\rvert^4=\E\left\lvert~ \E_{k} a_{k1}-\E_k d_k~\right\rvert^4 \leq ~C_0\lb \E ~g_k^4+\E ~h_k^4\rb= O(n^{-2}),
\end{align*}
which implies that for any $\veps>0$,
\begin{align*}
&~\sum_{k=1}^p\E\lb \Big(\sum_{i=1}^r\alpha_i Y_k(z_i)\Big)^2\mathds{1}_{\left\{\left\lvert\sum_{i=1}^r\alpha_i Y_k(z_i)\right\rvert\geq \veps\right\}}\rb\leq \frac{1}{\veps^2}\sum_{k=1}^{p}\E \Big\lvert \sum_{i=1}^r\alpha_i Y_k(z_i)\Big\rvert^4\\
\leq &~\frac{p~C_0}{\veps^2}~  \E\left\lvert ~\E_{k}\frac{\d~}{\d z_i}\lb \widetilde{b}_k(z_i)g_k(z_i)\rb~\right\rvert^4 =O(n^{-1}),
\end{align*}
then condition \eqref{eq:MCLTC1} is verified.

Therefore, the remaining is devoted to find the limit of
\eqref{limito}.
Denote $\D_k^{-1}(z)=\lb \mathbf{\Thta}^{\T}_k\mathbf{\Thta}_k-z\I_M\rb^{-1}$, then $\B_k=\I_M+z\D_k^{-1}(z)$, so we have
\begin{align*}
\E\left\lvert \widetilde{b}_k(z)-b_n(z)\right\rvert^2
=\E ~\widetilde{b}_k^2(z)b_n^2(z)\left[\frac{1}{3M}\tr\lb \TT\B_k\TT^{\T}+\B_k\rb-\frac{1}{3M}\E\tr\lb \TT\B_k\TT^{\T}+\B_k\rb\right]^2=O(n^{-2}),
\end{align*}
which implies
\begin{align*}
\sum_{k=1}^{p}\E_{k-1}\left[\E_k\lb \widetilde{b}_k(z_1)g_k(z_1)\rb \E_k\lb \widetilde{b}_k(z_2) g_k(z_2)\rb\right]-b_n(z_1)b_n(z_2)\sum_{k=1}^p\E_{k-1}\left[\E_kg_k(z_1)~\E_k g_k(z_2)\right]\xrightarrow{p}0.
\end{align*}
Therefore, it remains to find the limit of 
\begin{align}\label{f1}
\frac{d}{dz_1}\frac{d}{dz_2}\left\{b_n(z_1)b_n(z_2)\sum_{k=1}^p\E_{k-1}\left[\E_kg_k(z_1)~\E_k g_k(z_2)\right]\right\},
\end{align}
which in turn gives  the limit of \eqref{limito}.

In fact, denote $\H_k(z)=\I_M-\E_k \B_k(z)=-z~ \E_k \D^{-1}_k(z) \in \sigma(\mathcal F_{k-1})$, we can write
\begin{align*}
&\quad ~\E_{k-1} \Big[ \E_k g_k(z_1) \E_k g_k(z_2)\Big]
=\E \Big[ (\vv^{\T}_k\H_k(z_1)\vv_k-\E\vv^{\T}_k\H_k(z_1)\vv_k)(\vv^{\T}_k\H_k(z_2)\vv_k-\E\vv^{\T}_k\H_k(z_2)\vv_k)\Big]~.
\end{align*}
By applying Lemma~\ref{lem:quadratics}, we have
\begin{align}
&\quad b_n(z_1)b_n(z_2)\sum_{k=1}^p\E_{k-1}\left[\E_kg_k(z_1)~\E_k g_k(z_2)\right]\nonumber\\
&=b_n(z_1)b_n(z_2)\sum_{k=1}^p\frac{2}{9M^2}\tr\lb \TT\H_k(z_1)\TT^{\T}\TT\H_k(z_2)\TT^{\T}\rb \label{d1}\\
&\quad -b_n(z_1)b_n(z_2)\sum_{k=1}^p\frac{2}{15M^2}\sum_{i=1}^n\lb\TT\H_k(z_1)\TT^{\T}\rb_{ii}\lb\TT\H_k(z_2)\TT^{\T}\rb_{ii}\label{d2} \\
&\quad +b_n(z_1)b_n(z_2)\sum_{k=1}^p\frac{4}{45M^2}\Big[\sum_{i<\ell}\lb\TT\H_k(z_1)\TT^{\T}\rb_{ii}\lb\TT\H_k(z_2)\rb_{\ell,(i\ell)}-\sum_{\ell<i}\lb\TT\H_k(z_1)\TT^{\T}\rb_{ii}\lb\TT\H_k(z_2)\rb_{\ell,(\ell i)}\Big] \label{d3}\\
&\quad +b_n(z_1)b_n(z_2)\sum_{k=1}^p\frac{4}{45M^2}\Big[\sum_{i<\ell}\lb\TT\H_k(z_2)\TT^{\T}\rb_{ii}\lb\TT\H_k(z_1)\rb_{\ell,(i\ell)}-\sum_{\ell<i}\lb\TT\H_k(z_2)\TT^{\T}\rb_{ii}\lb\TT\H_k(z_1)\rb_{\ell,(\ell i)}\Big] \label{d4}\\
&\quad +b_n(z_1)b_n(z_2)\sum_{k=1}^p\frac{4}{45M^2}\Big[\sum_{j<i<t}\lb\TT\H_k(z_1)\rb_{i,(it)}\lb\TT\H_k(z_2)\rb_{j,(jt)}+\sum_{t<j<i}\lb\TT\H_k(z_1)\rb_{i,(ti)}\lb\TT\H_k(z_2)\rb_{j,(tj)}\Big] \label{d5}\\
&\quad +b_n(z_1)b_n(z_2)\sum_{k=1}^p\frac{4}{45M^2}\Big[\sum_{s<i<j}\lb\TT\H_k(z_1)\rb_{i,(si)}\lb\TT\H_k(z_2)\rb_{j,(sj)}  +\sum_{i<j<s}\lb\TT\H_k(z_1)\rb_{i,(is)}\lb\TT\H_k(z_2)\rb_{j,(js)}\Big] \label{d6}\\
&\quad -b_n(z_1)b_n(z_2)\sum_{k=1}^p\frac{4}{45M^2}\Big[\sum_{j<t<i}\lb\TT\H_k(z_1)\rb_{i,(ti)}\lb\TT\H_k(z_2)\rb_{j,(jt)}+\sum_{i<s<j}\lb\TT\H_k(z_1)\rb_{i,(is)}\lb\TT\H_k(z_2)\rb_{j,(sj)}\Big] \label{d7}\\
&\quad +O(n^{-1})~.\nonumber
\end{align}
Therefore, in order to obtain the limiting variance of $M^1_n(z)$, we need to find out the limit of each terms in the expansion of 
$b_n(z_1)b_n(z_2)\sum_{k=1}^p\E_{k-1}\left[\E_kg_k(z_1)~\E_k g_k(z_2)\right]$, i.e, the limit of \eqref{d1}$\sim$\eqref{d7}.
 The main technical point of derivation for the limit of these terms is to replace $\H_k(z)$ with $z\lb -\widehat{\TT}_N(z)\rb^{-1}$ and figure out the orders of the remainder. Specifically, by making use of the results in Lemma~\ref{lem:generalvkbound}$\sim$\ref{lem:quadratics}, we have the limit for \eqref{d1}$\sim$\eqref{d7} given in the following lemma, whose proof  is relegated to the supplementary file. 

\begin{lemma}\label{lem:fan1}
	With the same notations as in the previous context,  \eqref{d1} can be approximated by 
\[
	\frac{2}{p}\sum_{k=1}^p\frac{A(z_1,z_2)}{1-\frac{k-1}{p}A(z_1,z_2)}+O\lb n^{-1/2}\rb,
\]
	where
	\begin{align*}
	A(z_1,z_2)&=\frac{m_{F^{c_n}}(z_1)m_{F^{c_n}}(z_2)pn^3}{9M^2 \lb 1+\frac{2}{3}c_n m_{F^{c_n}}(z_1) \rb\lb 1+\frac{2}{3}c_n m_{F^{c_n}}(z_2)\rb}
	=\frac{4c_nm_{F^{c_n}}(z_1)m_{F^{c_n}}(z_2)}{m_{F^{c_n}}(z_1)-m_{F^{c_n}}(z_2)};
	\end{align*}
	\begin{align*}
	\eqref{d2}&=-\frac{8c_nm_{F^{c_n}}(z_1)m_{F^{c_n}}(z_2)}{15\lb 1+\frac{2}{3}c_nm_{F^{c_n}}(z_1) \rb \lb 1+\frac{2}{3}c_nm_{F^{c_n}}(z_2) \rb}+o_p(1);
	\end{align*}
\eqref{d3} and \eqref{d4} share the same limit, which is 
\[-\frac{16c_nm_{F^{c_n}}(z_1)m_{F^{c_n}}(z_2)}{45\lb 1+\frac{2}{3}c_nm_{F^{c_n}}(z_1)\rb\lb 1+\frac{2}{3}c_nm_{F^{c_n}}(z_2)\rb};\]
all the remaining  terms \eqref{d5}, \eqref{d6} and \eqref{d7}
share the same limit, which is
\[\frac{16c_nm_{F^{c_n}}(z_1)m_{F^{c_n}}(z_2)}{135\lb 1+\frac{2}{3}c_nm_{F^{c_n}}(z_1)\rb\lb 1+\frac{2}{3}c_nm_{F^{c_n}}(z_2)\rb}.\]
\end{lemma}

Combining the limits of  \eqref{d1}$\sim$\eqref{d7}  using Lemma~\ref{lem:fan1} and then taking derivatives with respect to $z_1$ and $z_2$ gives the limit of \eqref{f1}, which is exactly given by \eqref{eq:lemVar}. To this end,
we have proved that  $(M^1_n(z_1)\cdots M^1_n(z_r))^{\T}$ is jointly Gaussian with covariance function $\cov(M^1_n(z_1),M^1_n(z_2))$ given by \eqref{eq:lemVar}.

\subsubsection{Step 2: Tightness of  $M_n^1(z)$}
To prove the tightness of $M_n^1(z)$, it is sufficient to prove the moment condition (12.51) of \cite{billingsley1968convergence}, i.e.
\begin{align}\label{summ}
\sup_{n; z_1, z_2 \in \mathcal{C}_n}\frac{\E~ \big|M_n^1(z_1)-M_n^1(z_2)\big|^2}{|z_1-z_2|^2}
\end{align}
is finite.
According to edge university of the Kendall's rank correlation matrix derived in \cite{BZG18}, we can assume $\D^{-1}(z)$, $\D^{-1}_i(z)$ and $\D^{-1}_{ij}(z)$ are all bounded in $n$ and $z \in \mathcal{C}_n$. Then by applying  Lemma \ref{lem:vkbound}, the verification of  \eqref{summ} follows similar  procedure as developed in \cite{BS04} and the details will be omitted here.

\subsubsection{Step 3:  Convergence of $M_n^2(z)$}
The proof of Lemma \ref{lem:main} will be complete if we could show that $\{M^2_n(z)\}$ is bounded and form an equicontinuous family and converges to some fixed constant, which is given by \eqref{eq:lemMean}. As for the  boundedness and equicontinuity, it is easy to verify following the steps in \cite{BS04} so the main task in this part is to derive the limit of $\{M^2_n(z)\}$.

Notice that for $m_{F^{c_n}}(z)$, it satisfies
\begin{gather*}
    \frac{1}{1+\frac23 c_n m_{F^{c_n}}(z)}=1-c_n-c_n(z-\frac{1}{3})m_{F^{c_n}}(z),
\end{gather*}
so we define
\begin{align*}
A_n&=\frac{1}{1+\frac23 c_n\E m_{F^{\K_n}}(z)}-1+c_n+c_n(z-\frac{1}{3})\E m_{F^{\K_n}}(z)~.
\end{align*}
Consider the following
\begin{align*}
m_{F^{c_n}}(z) A_n
&=c_n\lb m_{F^{c_n}}(z)-\E m_{F^{\K_n}}(z)\rb +c_n \E m_{F^{\K_n}}(z)\left[ 1+(z-\frac{1}{3})m_{F^{c_n}}(z) \right]\\
&\quad -\frac{2c_nm_{F^{c_n}}(z)\cdot\E m_{F^{\K_n}}(z) }{3+2c_n\E m_{F^{\K_n}}(z)}\\
&=c_n\lb m_{F^{c_n}}(z)-\E m_{F^{\K_n}}(z)\rb -  \frac{4c^2_nm_{F^{c_n}}(z)\E m_{F^{\K_n}}(z)\lb m_{F^{c_n}}(z)-\E m_{F^{\K_n}}(z)\rb }{\lb 3+2c_n\E m_{F^{\K_n}}(z)\rb\lb 3+2c_n m_{F^{c_n}}(z)\rb },
\end{align*}
so we have
\begin{align}\label{emm}
\E m_{F^{\K_n}}(z)- m_{F^{c_n}}(z) &=m_{F^{c_n}}(z) ~A_n\left[-c_n+\frac{4c_n^2~m_{F^{c_n}}(z)\cdot\E m_{F^{\K_n}}(z)}{\lb 3+2 c_n\E m_{F^{\K_n}}(z)\rb\lb 3+ 2 c_n m_{F^{c_n}}(z)\rb}\right]^{-1}\nonumber\\
&=\cfrac{m_{F^{c_n}}(z) \left[ 1+(z-\frac13)\E m_{F^{\K_n}}(z)-\frac{2 \E m_{F^{\K_n}}(z) }{3+2c_n\E m_{F^{\K_n}}(z)}\right]\quad}{-1+\frac{4c_nm_{F^{c_n}}(z)~\E m_{F^{\K_n}}(z)}{\big( 3+2 c_n\E m_{F^{\K_n}}(z)\big)\big( 3+ 2 c_n m_{F^{c_n}}(z)\big)}}.
\end{align}

On the other hand, with some non-trivial calculations, 
\begin{align}\label{trdiff}
&\quad \tr\lb\E \D^{-1}(z)\rb-\tr\lb -\widehat{\TT}_N(z)\rb^{-1}\\
&=p \E m_{F^{\K_n}}(z)+\frac{p-M}{z}+\frac{1}{z}\tr\lb \I_M+\frac{p-1}{3M} \E m_{F^{\K_n}}(z)\lb \TT^{\T}\TT+\I_M\rb\rb^{-1}\nonumber\\
&=\frac{p}{z} \left[ 1+(z-\frac13)\E m_{F^{\K_n}}(z)-\frac{2 \E m_{F^{\K_n}}(z) }{3+2c_n\E m_{F^{\K_n}}(z)}\right]+D_n(z)+o(1)~,\nonumber
\end{align}
where $D_n(z)$ is a notation to denote the constant part, i.e.
\begin{align}\label{dn}
D_n(z)=\frac{2c_n^2\lb\E m_{F^{\K_n}}(z)\rb^2}{9z}+\frac{\E m_{F^{\K_n}}(z)}{3z}+\frac{\frac{8}{9}c_n^3\lb\E m_{F^{\K_n}}(z)\rb^3+4c_n^2\lb\E m_{F^{\K_n}}(z)\rb^2+2\E m_{F^{\K_n}}(z)}{3z\lb 1+\frac{2}{3}c_n\E m_{F^{\K_n}}(z)\rb^2}~.
\end{align}

 Therefore, combining \eqref{emm} and \eqref{trdiff}, we have
\begin{align}\label{eq:diff}
&p\lb\E m_{F^{\K_n}}(z)- m_{F^{c_n}}(z)\rb= \frac{z m_{F^{c_n}}(z)\lb \tr\lb\E \D^{-1}(z)\rb-\tr\lb -\widehat{\TT}_N(z)\rb^{-1}-D_n(z)\rb}{-1+\cfrac{\frac{4}{9}c_n~m_{F^{c_n}}(z)\cdot\E m_{F^{\K_n}}(z)}{\lb 1+\frac{2}{3} c_n\E m_{F^{\K_n}}(z)\rb\lb 1+ \frac{2}{3}c_n m_{F^{c_n}}(z)\rb}}~+o(1)~.\\ \nonumber
\end{align}
Thus, what remains is to find the limit of the term on the right hand side of Equation \eqref{eq:diff}. First, considering the part involving $\tr\lb -\widehat{\TT}_N(z)\rb^{-1}- \tr\lb\E \D^{-1}(z)\rb$, we have
\begin{align*}
&\tr\lb -\widehat{\TT}_N(z)\rb^{-1}- \tr\lb\E \D^{-1}(z)\rb\\
=&~\tr \left[ \lb -\widehat{\TT}_N(z)\rb^{-1}\E\lb \sum_{k=1}^p\vv_k\vv_k^{\T}+\frac{p-1}{3M}z\E  m_{F^{\K_n}}(z)\lb \TT^{\T}\TT+\I_M\rb\rb\D^{-1}(z)\right]\\
=&~\sum_{k=1}^p\E z\beta_k\vv_k^{\T}\D_k^{-1}(z)\lb -\frac{p-1}{3M}z\E m_{F^{\K_n}}(z)\lb \TT^{\T}\TT+\I_M\rb-z\I_M\rb^{-1}\vv_k\\
&-\sum_{k=1}^p z\E\beta_k\cdot\tr\left[ \frac{1}{3M}\lb \TT^{\T}\TT+\I_M\rb\E\D^{-1}(z)\lb -\frac{p-1}{3M}z\E m_{F^{\K_n}}(z)\lb \TT^{\T}\TT+\I_M\rb-z\I_M\rb^{-1}\right]\\
&+ z\E  \beta_k \cdot\tr\left[ \frac{1}{3M}\lb \TT^{\T}\TT+\I_M\rb\E\D^{-1}(z)\lb -\widehat{\TT}_N(z)\rb^{-1}\right]\\
=&~\sum_{k=1}^{p} \E z\beta_k\left[ \vv_k^{\T}\D_k^{-1}(z)\lb -\widehat{\TT}_N(z)\rb^{-1}\vv_k-\frac{1}{3M}\E\tr \D^{-1}(z)\lb -\widehat{\TT}_N(z)\rb^{-1}\lb \TT^{\T}\TT+\I_M\rb\right]\\
&- z\E  m_{F^{\K_n}}(z) \cdot\tr\left[ \frac{1}{3M}\lb \TT^{\T}\TT+\I_M\rb\E\D^{-1}(z)\lb -\widehat{\TT}_N(z)\rb^{-1}\right].
\end{align*}
Define $r_k=z\left[ \vv_k^{\T}\D_k^{-1}(z)\vv_k-\frac{1}{3M}\E \tr \D_k^{-1}(z)\lb \TT^{\T}\TT+\I_M\rb\right]$, we have
\[
\beta_k=b_n-b_n\beta_kr_k=b_n-b_n^2r_k+\beta_kb_n^2r_k^2.
\]
Then
\begin{align*}
&\frac{1}{3M}\E\tr\lb\D^{-1}(z)-\D_k^{-1}(z)\rb\lb -\widehat{\TT}_N(z)\rb^{-1}\lb \TT^{\T}\TT+\I_M\rb\\
=&-\frac{1}{3M}\E z\beta_k\vv_k^{\T}\D_k^{-1}(z)\lb -\widehat{\TT}_N(z)\rb^{-1}\lb \TT^{\T}\TT+\I_M\rb\D_k^{-1}(z)\vv_k\\
=&-zb_n\E\tr\left[ \D_k^{-1}(z) \lb -\widehat{\TT}_N(z)\rb^{-1}\frac{1}{3M}\lb \TT^{\T}\TT+\I_M\rb\D_k^{-1}(z)\frac{1}{3M}\lb \TT^{\T}\TT+\I_M\rb\right]+o(\frac{1}{n}),
\end{align*}
and
\begin{align}\label{diff1}
&\tr\lb -\widehat{\TT}_N(z)\rb^{-1}- \tr\lb\E \D^{-1}(z)\rb\\
=&~\sum_{k=1}^{p} \E z\beta_k\left[ \vv_k^{\T}\D_k^{-1}(z)\lb -\widehat{\TT}_N(z)\rb^{-1}\vv_k-\frac{1}{3M}\E\tr \D_k^{-1}(z)\lb -\widehat{\TT}_N(z)\rb^{-1}\lb \TT^{\T}\TT+\I_M\rb\right]\nonumber\\
&- zb_n^2\sum_{k=1}^p\E\tr\left[ \D_k^{-1}(z) \lb\I_M+\frac{p-1}{3M}\E m_{F^{\K_n}}(z)\lb \TT^{\T}\TT+\I_M\rb\rb^{-1}\frac{1}{3M}\lb \TT^{\T}\TT+\I_M\rb\D_k^{-1}(z)\frac{1}{3M}\lb \TT^{\T}\TT+\I_M\rb\right]\nonumber\\ \nonumber
&- z\E  m_{F^{\K_n}}(z) \cdot\tr\left[ \frac{1}{3M}\lb \TT^{\T}\TT+\I_M\rb\E\D^{-1}(z)\lb -\widehat{\TT}_N(z)\rb^{-1}\right]. 
\end{align}
As for the first term in \eqref{diff1}, we have
\begin{align*}
&\sum_{k=1}^{p} \E z\beta_k\left[ \vv_k^{\T}\D_k^{-1}(z)\lb -\widehat{\TT}_N(z)\rb^{-1}\vv_k-\frac{1}{3M}\E\tr \D_k^{-1}(z)\lb -\widehat{\TT}_N(z)\rb^{-1}\lb \TT^{\T}\TT+\I_M\rb\right]\\
=&-b_n^2\sum_{k=1}^{p}\E zr_k\vv_k^{\T}\D_k^{-1}(z)\lb -\widehat{\TT}_N(z)\rb^{-1}\vv_k\\
&~+b_n^2\sum_{k=1}^{p}\E z\beta_kr_k^2\left[ \vv_k^{\T}\D_k^{-1}(z)\lb -\widehat{\TT}_N(z)\rb^{-1}\vv_k-\frac{1}{3M}\E\tr \D_k^{-1}(z)\lb -\widehat{\TT}_N(z)\rb^{-1}\lb \TT^{\T}\TT+\I_M\rb \right]\\
=&-b_n^2\sum_{k=1}^{p}\E zr_k\vv_k^{\T}\D_k^{-1}(z)\lb -\widehat{\TT}_N(z)\rb^{-1}\vv_k+b_n^2\sum_{k=1}^{p} \left\{\E z\beta_kr_k^2\vv_k^{\T}\D_k^{-1}(z)\lb -\widehat{\TT}_N(z)\rb^{-1}\vv_k\right.\\
&\left. -\E z\beta_kr_k^2\tr \lb \D_k^{-1}(z)\lb -\widehat{\TT}_N(z)\rb^{-1}\frac{1}{3M}\lb \TT^{\T}\TT+\I_M\rb \rb+\cov \lb z\beta_kr_k^2,\tr \left[ \D_k^{-1}(z)\lb -\widehat{\TT}_N(z)\rb^{-1}\frac{1}{3M}\lb \TT^{\T}\TT+\I_M\rb \right]\rb\right\}\\
=& zb_n^2\sum_{k=1}^{p} \E\left[ \vv_k^{\T} \D_k^{-1}(z)\vv_k-\frac{1}{3M}\E \tr\D_k^{-1}(z)\lb \TT^{\T}\TT+\I_M\rb\right]\cdot\left[ \vv_k^{\T}\D_k^{-1}(z)\lb \I_M+\frac{p-1}{3M}\E m_{F^{\K_n}}(z)\lb \TT^{\T}\TT+\I_M\rb \rb^{-1}\vv_k\right.\\
&~\left. -\frac{1}{3M}\tr\lb \D_k^{-1}(z)\lb \I_M+\frac{p-1}{3M}\E m_{F^{\K_n}}(z)\lb \TT^{\T}\TT+\I_M\rb \rb^{-1}\lb \TT^{\T}\TT+\I_M\rb\rb\right]+o(1).
\end{align*}
Denote $\WW^{-1}(z)=\lb \I_M+\frac{p-1}{3M}\E m_{F^{\K_n}}(z)\lb \TT^{\T}\TT+\I_M\rb \rb^{-1}$ for short, then using Lemma~\ref{lem:quadratics}, we have
\begin{align}
&\quad \tr\lb -\widehat{\TT}_N(z)\rb^{-1}- \tr\lb\E \D^{-1}(z)\rb\nonumber\\
&= -z\E  m_{F^{\K_n}}(z)\tr\left[ \frac{1}{3M}\lb \TT^{\T}\TT+\I_M\rb\E\D^{-1}(z)\lb -\widehat{\TT}_N(z)\rb^{-1}\right] \label{td1}\\
+&zb_n^2 \frac{1}{9M^2} \sum_{k=1}^p \E \tr\lb \TT \D_k^{-1}(z)\TT^{\T}\TT \D_k^{-1}(z)\WW^{-1}(z)\TT^{\T}\rb \label{td2}\\
-&zb_n^2\frac{2}{15M^2}\sum_{k=1}^p  \sum_{i=1}^n \E\lb\TT\D_k^{-1}(z)\TT^{\T}\rb_{ii}\lb\TT\D_k^{-1}(z)\WW^{-1}(z)\TT^{\T}\rb_{ii}  \label{td3}\\
+&~zb_n^2\frac{4}{45M^2}\sum_{k=1}^p \sum_{i<\ell} \E\left[\lb \TT\D_k^{-1}(z)\TT^{\T}\rb_{ii}\lb \TT\D_k^{-1}(z)\WW^{-1}(z)\rb_{\ell,(i\ell)}-\lb \TT\D_k^{-1}(z)\WW^{-1}(z)\TT^{\T}\rb_{\ell \ell}\lb \TT\D_k^{-1}(z)\rb_{i,(i \ell)}\right]\label{td4}\\
+&~zb_n^2\frac{4}{45M^2}\sum_{k=1}^p \sum_{i<\ell} \E \left[\lb \TT\D_k^{-1}(z)\WW^{-1}(z)\TT^{\T}\rb_{ii}\lb \TT\D_k^{-1}(z)\rb_{\ell,(i \ell)}-\lb \TT\D_k^{-1}(z)\TT^{\T}\rb_{\ell \ell}\lb \TT\D_k^{-1}(z)\WW^{-1}(z)\rb_{i,(i \ell)}\right]\nonumber\\
+&~zb_n^2\sum_{k=1}^p\frac{4}{45M^2}\E\left[\sum_{j<i<t}\lb\TT\D_k^{-1}(z)\rb_{i,(it)}\lb\TT\D_k^{-1}(z)\WW^{-1}(z)\rb_{j,(jt)}+\sum_{j<i<t}\lb\TT\D_k^{-1}(z)\WW^{-1}(z)\rb_{i,(it)}\lb\TT\D_k^{-1}(z)\rb_{j,(jt)}\right] \label{td5}\\
+&~zb_n^2\sum_{k=1}^p\frac{4}{45M^2}\E\left[\sum_{t<j<i}\lb\TT\D_k^{-1}(z)\rb_{i,(ti)}\lb\TT\D_k^{-1}(z)\WW^{-1}(z)\rb_{j,(tj)}+\sum_{t<j<i}\lb\TT\D_k^{-1}(z)\WW^{-1}(z)\rb_{i,(ti)}\lb\TT\D_k^{-1}(z)\rb_{j,(tj)}\right] \nonumber\\
-&~zb_n^2\sum_{k=1}^p\frac{4}{45M^2}\E\left[\sum_{j<t<i}\lb\TT\D_k^{-1}(z)\rb_{i,(ti)}\lb\TT\D_k^{-1}(z)\WW^{-1}(z)\rb_{j,(jt)}+\sum_{j<t<i}\lb\TT\D_k^{-1}(z)\WW^{-1}(z)\rb_{i,(ti)}\lb\TT\D_k^{-1}(z)\rb_{j,(jt)}\right]\nonumber\\
~+&o(1). \nonumber
\end{align}

Next, we look into each terms above  by replacing $\D^{-1}(z)$ with $\lb -\widehat{\TT}_N(z)\rb^{-1}$ and figure out the orders of the remainder. Through some non-trivial calculations, we have the following lemma, whose proof is relegated to the supplement file.

\begin{lemma}\label{lem:fan2}
With the same notations as in the previous context, we have
\begin{gather*}
\eqref{td1}=-\frac{\E  m_{F^{\K_n}}(z)}{3z}\lb 1+\frac{2}{\lb 1+\frac{2}{3}c_n \E  m_{F^{\K_n}}(z)\rb^2}\rb+o(1),\\
\eqref{td2}
  =\frac{4c_n b^2_n}{9z\lb 1+\frac{2}{3}c_n \E m_{F^{\K_n}}(z)\rb^3}\left[1-\frac{4c_n\lb \E m_{F^{\K_n}}(z)\rb^2}{9\lb 1+\frac{2}{3}c_n \E m_{F^{\K_n}}(z)\rb^2}\right]^{-1}+o(1),\\
\eqref{td3}=-\frac{8c_n b^2_n}{15z\lb 1+\frac{2}{3}c_n\E  m_{F^{\K_n}}(z)\rb^3}+o(1).
  \end{gather*} 
  Moreover, 
the four terms in \eqref{td4} share the same limit, which is
  $$-\frac{8c_n b^2_n}{45z\lb 1+\frac{2}{3}c_n\E m_{F^{\K_n}}(z)\rb^3}$$
and all the six terms in \eqref{td5}
share the same limit, which is
 $$\frac{8 c_n b^2_n}{135z\lb 1+\frac{2}{3}c_n\E m_{F^{\K_n}}(z)\rb^3}.$$
\end{lemma}

Collecting all the limits in  Lemma~\ref{lem:fan2}, we have
\begin{align}\label{tod}
&\tr\lb -\widehat{\TT}_N(z)\rb^{-1}- \tr\lb\E \D^{-1}(z)\rb=-\frac{\E m_{F^{\K_n}}(z)}{3z}\lb 1+\frac{2}{\lb 1+\frac{2}{3}c_n \E m_{F^{\K_n}}(z)\rb^2}\rb\\
+&\frac{4c_n\lb \E m_{F^{\K_n}}(z)\rb^2}{9z\lb 1+\frac{2}{3}c_n \E m_{F^{\K_n}}(z)\rb^3}\left[1-\frac{4c_n\lb \E m_{F^{\K_n}}(z)\rb^2}{9\lb 1+\frac{2}{3}c_n \E m_{F^{\K_n}}(z)\rb^2}\right]^{-1}-\frac{8c_n\lb \E m_{F^{\K_n}}(z)\rb^2}{9z \lb 1+\frac{2}{3}c_n\E  m_{F^{\K_n}}(z)\rb^3}+o(1)~. \nonumber
\end{align}

Finally, combing \eqref{dn}, \eqref{eq:diff} and \eqref{tod},  the limit of $\{M^2_n(z)\}$  is found to be given by \eqref{eq:lemMean}.

 To this end, the proof of Lemma  \ref{lem:main} is complete.

\medskip

\section{Some useful lemmas}

\begin{definition}
    Let $X^{(n)}$ and $Y^{(n)}$ be two sequences of nonnegative random variables, we say that  $Y^{(n)}$ \em{stochastically dominates} $X^{(n)}$ if for all small $\veps>0$ and large $d>0$,
    \[\mathbb{P}\lb X^{(n)}>n^{\veps}Y^{(n)}\rb\leq n^{-d},\]
    for sufficiently large $n\geq n_0(\veps,d)$ and we denote as $X^{(n)}\prec Y^{(n)}$.
\end{definition}

\begin{lemma}\label{lem:bili08}
$[\mbox{Theorem 35.12 in \cite{bili08}}].$  Suppose for each $n$, $$Y_{n1},Y_{n2},\cdots,Y_{nr_n}$$ is a real martingale difference sequence with respect to the increasing $\sigma-$field $\{\mathcal{F}_{nj}\}$ having second moments. If as $n\rightarrow \infty$,
    \[ \sum_{j=1}^{r_n}\E\lb Y_{nj}^2~|~\mathcal{F}_{n,j-1}\rb\xrightarrow{p}\sigma^2,\]
    where $\sigma^2$ is a positive constant and for each $\veps>0$,
    \[\sum_{j=1}^{r_n}\E\lb Y_{nj}^2~\mathds{1}_{\{|Y_{nj}|\geq\veps\}}\rb\xrightarrow{p}0,\]
    then
    \[\sum_{j=1}^{r_n} Y_{nj}\xrightarrow{d} N(0,\sigma^2).\]
\end{lemma}

\begin{lemma}\label{lem:generalvkbound}
For $\vv_k$ as defined in \eqref{vk} and any $M\times M$ deterministic Hermitian matrices  $\A_k$, $\B_{\ell}$, $1\leq k\leq m$, $1\leq\ell\leq q$, under the same assumptions as in Theorem \ref{thm:mainclt}, there exist constant $K>0$ such that
\begin{equation}\label{eq:vkgeneralbound}
 \E\left\{ \prod_{k=1}^{m}\vv_t^{\T}\A_k\vv_t\cdot\prod_{\ell=1}^{q}\left[\vv_t^{\T}\B_{\ell}\vv_t-\tr\lb \B_{\ell}\lb \frac{1}{3M}\TT^{\T}\TT+\frac{1}{3M}\I_M\rb\rb\right]\right\}\leq Kn^{-\frac{q}{2}}\prod_{k=1}^{m}||\A_k||\prod_{\ell=1}^{q}||\B_{\ell}||,
 \end{equation}
where $||\cdot ||$ denotes the spectral norm.
\end{lemma}

\begin{lemma}\label{lem:vkbound}
For the $M$-dimensional random vector $\vv_k$ as defined in \eqref{vk} and any  $M\times M$ deterministic Hermitian matrix $\A$ which has bounded spectral norm, under the same assumptions as in Theorem \ref{thm:mainclt}, we have
\[\E\left\lvert \vv_k^{\T}\A\vv_k-\frac{1}{3M}\tr\lb \TT\A\TT^{\T}\rb-\frac{1}{3M}\tr\A\right\rvert^m= O\lb n^{-\frac{m}{2}}\rb.\]
\end{lemma}

\medskip


\begin{lemma}\label{lem:quadratics}
    Let $\A,~\B$ be any $M\times M$ deterministic symmetric matrices, we have
    \begin{align*}
    &~\E \Big(\vv_k^{\T}\A\vv_k-\E\vv_k^{\T}\A\vv_k\Big)\Big(\vv_k^{\T}\B\vv_k-\E \vv_k^{\T}\B\vv_k\Big)\\
    =&~\E\lb \vv_k^{\T}\A\vv_k\vv_k^{\T}\B\vv_k\rb-\lb\frac{1}{3M}\tr\lb \TT\A\TT^{\T}\rb+\frac{1}{3M}\tr(\A)\rb\lb\frac{1}{3M}\tr\lb \TT\B\TT^{\T}\rb+\frac{1}{3M}\tr(\B)\rb\\
    =&~\frac{2}{9M^2}\tr\lb \TT\A\TT^{\T}\TT\B\TT^{\T}+\A\B\rb-\frac{2}{15M^2}\sum_{i=1}^n\left[\lb\TT\A\TT^{\T}\rb_{ii}\lb\TT\B\TT^{\T}\rb_{ii}+\A_{ii}\B_{ii}\right]+ \frac{4}{9M^2}\tr\lb \TT\A\B\TT^{\T}\rb\\
    +&~\frac{4}{45M^2}\sum_{i<\ell}\left[\lb\TT\A\TT^{\T}\rb_{ii}\lb\TT\B\rb_{\ell,(i\ell)}+\lb\TT\B\TT^{\T}\rb_{ii}\lb\TT\A\rb_{\ell,(i\ell)}\right]-\frac{4}{45M^2}\sum_{\ell<i}\left[\lb\TT\A\TT^{\T}\rb_{ii}\lb\TT\B\rb_{\ell,(\ell i)}+\lb\TT\B\TT^{\T}\rb_{ii}\lb\TT\A\rb_{\ell,(\ell i)}\right]\\
        +&~\frac{8}{45M^2}\sum_{i<j}\left[\lb\TT\A\TT^{\T}\rb_{ij}\left[ \lb\TT\B\rb_{i,(ij)}-\lb\TT\B\rb_{j,(ij)}\right]+\lb\TT\B\TT^{\T}\rb_{ij}\left[ \lb\TT\A\rb_{i,(ij)}-\lb\TT\A\rb_{j,(ij)}\right]\right]\\
    +&~\frac{4}{45M^2}\left[\sum_{j<i<t}\lb\TT\A\rb_{i,(it)}\lb\TT\B\rb_{j,(jt)}+\sum_{t<j<i}\lb\TT\A\rb_{i,(ti)}\lb\TT\B\rb_{j,(tj)}-\sum_{j<t<i}\lb\TT\A\rb_{i,(ti)}\lb\TT\B\rb_{j,(jt)}\right]\\
        +&~\frac{4}{45M^2}\left[\sum_{j<i<t}\lb\TT\B\rb_{i,(it)}\lb\TT\A\rb_{j,(jt)}+\sum_{t<j<i}\lb\TT\B\rb_{i,(ti)}\lb\TT\A\rb_{j,(tj)}-\sum_{j<t<i}\lb\TT\B\rb_{i,(ti)}\lb\TT\A\rb_{j,(jt)}\right]\\
    +&~\frac{8}{45M^2}\sum_{1\leq i<j\leq n}\left[\lb\TT\A\rb_{i,(ij)}\lb\TT\B\rb_{j,(ij)}+ \lb\TT\B\rb_{i,(ij)}\lb\TT\A\rb_{j,(ij)}\right]\\
    +&~\frac{4}{45M^2}\left[\sum_{t<i<j}\lb\TT\A\rb_{i,(tj)}\lb\TT\B\rb_{j,(ti)}+\sum_{i<j<t}\lb\TT\A\rb_{i,(jt)}\lb\TT\B\rb_{j,(it)}-\sum_{i<t<j}\lb\TT\A\rb_{i,(tj)}\lb\TT\B\rb_{j,(it)}\right]\\
        +&~\frac{4}{45M^2}\left[\sum_{t<i<j}\lb\TT\B\rb_{i,(tj)}\lb\TT\A\rb_{j,(ti)}+\sum_{i<j<t}\lb\TT\B\rb_{i,(jt)}\lb\TT\A\rb_{j,(it)}-\sum_{i<t<j}\lb\TT\B\rb_{i,(tj)}\lb\TT\A\rb_{j,(it)}\right]\\
    +&~\frac{4}{45M^2}\sum_{i<j} \left[ \lb \TT\A\TT^{\T}\rb_{ij}\B_{(ij),(ij)}+ \lb \TT\B\TT^{\T}\rb_{ij}\A_{(ij),(ij)}\right] -\frac{4}{45M^2}\sum_{i<u<j}\left[\lb \TT\A\TT^{\T}\rb_{ij}\B_{(uj),(iu)}+\lb \TT\B\TT^{\T}\rb_{ij}\A_{(uj),(iu)}\right]\\
        +&~\frac{4}{45M^2}\sum_{i<j<t}\left[\lb \TT\A\TT^{\T}\rb_{ij}\B_{(jt),(it)}+\lb \TT\B\TT^{\T}\rb_{ij}\A_{(jt),(it)}\right]
        +\frac{4}{45M^2}\sum_{s<i<j}\left[\lb\TT\A\TT^{\T}\rb_{ij}\B_{(si),(sj)}+\lb \TT\B\TT^{\T}\rb_{ij}\A_{(si),(sj)}\right].
    \end{align*}
\end{lemma}

%
%
%

\bibliographystyle{apalike}
\bibliography{KEND}

\end{document}